\documentclass[reqno]{amsart}

\usepackage{amsmath,amsthm,amssymb,amsfonts}
\newcommand{\sgn}{\operatorname{sgn}}
\newcommand{\pfaff}{\mathop{\mathrm{pfaff}}}
\newcommand{\ep}{\varepsilon}
\numberwithin{equation}{section} 
\theoremstyle{plain}
\newtheorem{theo+}           {Theorem}      [section]
\newtheorem{prop+}  [theo+]  {Proposition}
\newtheorem{coro+}  [theo+]  {Corollary}
\newtheorem{lemm+}  [theo+]  {Lemma}
\newtheorem{defi+}  [theo+]  {Definition}

\theoremstyle{definition}
\newtheorem{exam+}  [theo+]  {Example}
\newtheorem{rema+}  [theo+]  {Remark}
\newtheorem{prob+}  [theo+]  {Problem}

\newenvironment{theorem}{\begin{theo+}}{\end{theo+}}
\newenvironment{proposition}{\begin{prop+}}{\end{prop+}}
\newenvironment{corollary}{\begin{coro+}}{\end{coro+}}
\newenvironment{lemma}{\begin{lemm+}}{\end{lemm+}}

\newenvironment{remark}{\begin{rema+}}{\end{rema+}}

\newenvironment{problem}  {\begin{prob+}}{\end{prob+}}

\begin{document}

\baselineskip 18pt
\larger[2]
\title
[Schur $Q$-polynomials, hypergeometric series and tableaux]
{Schur $Q$-polynomials, multiple hypergeometric series\\
 and enumeration of  marked shifted tableaux} 
\author{Hjalmar Rosengren}
\address
{Department of Mathematical Sciences
\\ Chalmers University of Technology and G\"oteborg
 University\\SE-412~96 G\"oteborg, Sweden}
\email{hjalmar@math.chalmers.se}
\urladdr{http://www.math.chalmers.se/{\textasciitilde}hjalmar}

\thanks{Research  supported by the Swedish Science Research
Council (Vetenskapsr\aa det)}

\begin{abstract}
We study Schur $Q$-polynomials evaluated on a geometric progression, or equivalently $q$-enumeration of marked shifted tableaux, seeking explicit formulas that remain regular at $q=1$. We obtain several such expressions as multiple basic hypergeometric series, and as determinants and pfaffians of continuous $q$-ultraspherical or continuous $q$-Jacobi polynomials. As special cases, we obtain simple closed formulas for staircase-type partitions.

\end{abstract}

\maketitle

\section{Introduction}  

The Schur $Q$-polynomials  originate in the work of
Schur \cite{s} on projective representations of the symmetric
group. Nowadays, they are recognized as the case $t=-1$ of the 
 Hall--Littlewood polynomials, which are in turn a special case of the
 Macdonald polynomials \cite{m}. 

In the present paper we  report on some investigations on Schur
$Q$-polynomials that were motivated by applications to sums
of squares \cite{rss}. In \cite{ro}, we used elliptic
pfaffian evaluations to derive, and generalize,
 certain triangular number identities
conjectured by Kac and Wakimoto \cite{kw} and proved first by Milne \cite{mi2,mi3,mi4}
and later by Zagier \cite{z}. Extending this work to the case of square numbers
leads naturally to quantities  related to Schur
$Q$-polynomials. To be precise, we 
note that, although one usually assumes that 
 the polynomial
$Q_\lambda$ is indexed by a partition $\lambda$, the   definition  \eqref{pdef}
makes sense  for general $\lambda\in\mathbb Z^m$. The quantities
alluded to may then be written $Q_{(\lambda,-\lambda)}(1,\dots,1)$,
where $\lambda$ is a partition. In  the present paper,  
however, we investigate $Q_\lambda(1^n)$ for general $\lambda$.

In many ways, Schur $Q$-polynomials and Schur polynomials have analogous properties. Therefore, it may be instructive to compare with the much simpler situation for the latter. The standard definition
$$s_\lambda(x)=\frac{\det_{1\leq i,j\leq n}(x_i^{\lambda_j+n-j})}{\det_{1\leq i,j\leq n}(x_i^{n-j})} $$
has a removable singularity at $x_1=\dots=x_n=1$. A well-known way to compute $s_\lambda(1^n)$ is by passing to the case when the variables are in geometric progression. Namely,  by 
 the Vandermonde determinant evaluation,
\begin{equation}\label{sfs}s_\lambda(1,q,\dots,q^{n-1})=
\frac{\prod_{1\leq i<j\leq n}(q^{\lambda_j+n-j}-q^{\lambda_i+n-i})}{\prod_{1\leq i<j\leq n}(q^{n-j}-q^{n-i})},
\end{equation}
which in the limit $q\rightarrow 1$ gives
$$s_\lambda(1^n)=\prod_{1\leq i<j\leq n}\frac{\lambda_i-\lambda_j+j-i}{j-i}.
$$

Our method for computing $Q_\lambda(1^n)$ is similar. In that case, $Q_\lambda(1,q,\dots,q^{n-1})$ is a multiple basic hypergeometric sum, see \eqref{phg}, again with a removable singularity at $q=1$. Since it does not simplify in general, the best one can hope for is to find a transformation formula, leading to a different sum without the apparent singularity. In the limit $q\rightarrow 1$ one would then obtain a formula for  $Q_\lambda(1^n)$ as a multiple classical hypergeometric sum. Finding such identities is our main goal.

It should be remarked that $Q_\lambda(1,q,\dots,q^{n-1})$ has a combinatorial meaning as a generating function (or $q$-enumeration) on marked shifted tableaux, see \eqref{qe}. Thus, we can equivalently formulate our goal as finding explicit solutions to the $q$-enumeration problem that remain regular in the limit $q\rightarrow 1$, corresponding to classical enumeration.

The key fact for finding transformation formulas with the desired property is 
Lemma \ref{tvp}, which allows us to identify the two-row $Q$-polynomial 
$$Q_{(\lambda_1,\lambda_2)}(1,q,\dots,q^{n})$$
 with the \emph{Christoffel--Darboux kernel} of certain continuous $q$-Jacobi polynomials. Like most of our results, it is formulated in terms of multivariable polynomials $P_n$ that are related to $Q_\lambda$ through
\begin{equation}\label{pq}Q_\lambda(1,q,\dots,q^n)=2^m P_n(q^{\lambda_1},\dots,q^{\lambda_m}).\end{equation}
The multi-row $Q$-polynomial similarly corresponds to a multivariable Christoffel--Darboux kernel. This fact  implies several pfaffian and determinantal formulas for the polynomials $P_n$. Besides the pfaffian formula \eqref{pfi}, which is essentially Schur's definition of the $Q$-polynomial, we have the transformed pfaffian formula of Corollary \ref{api} and the determinant formula of Corollary \ref{dt}.  In the limit $q\rightarrow 1$, the latter leads to an elegant determinant formula for $Q_\lambda(1^n)$, see Corollary \ref{fdc}.

We then use our new pfaffian and determinantal formulas to obtain hypergeometric 
formulas. 
Corollary \ref{api} allows us to express $P_n$ as a pfaffian of double sums. Series manipulation then leads to the multiple hypergeometric sum  of Theorem \ref{th}.
 Similarly, Corollary \ref{dt} expresses $P_n$ as a determinant of single sums, which leads to the Schlosser-type hypergeometric sums of Theorems \ref{ot} and
 \ref{dft}. In the limit $q\rightarrow 1$, Theorems \ref{th} and \ref{dft} yield  hypergeometric formulas for $Q_\lambda(1^n)$, thus fulfilling our main goal. (Theorem \ref{ot} retains the singularity at $q=1$.)

In some special cases, our expressions for  $Q_\lambda(1,q,\dots,q^{n-1})$ simplify, giving completely factored expressions similarly as for Schur polynomials.
\begin{itemize}
\item $n=\infty$, see Section \ref{ks} and Remark \ref{kr}. The explicit evaluation of 
$$Q_\lambda(1,q,q^2,\dots)$$ 
follows from an identity of Kawanaka \cite{k}, who does not, however, mention the relation to $Q$-polynomials.
\item $\lambda$ staircase, see Remark \ref{sr}. If $\lambda=(m,m-1,\dots,1)$, then $Q_\lambda=2^ms_\lambda$, so  $Q_\lambda(1,q,\dots,q^{n-1})$ factors by \eqref{sfs}.
\item $\lambda$ odd staircase, see Section \ref{dss}, Remark  \ref{oscr} and Remark \ref{sr}. If  $\lambda=(2m-1,2m-3,\dots,1)$, then $Q_\lambda(1,q,\dots,q^{n-1})$ can be evaluated using a discrete Selberg integral due to Milne \cite{mi}, or in the form needed here to Krattenthaler \cite{kr}. Again, the connection to $Q$-polynomials seems not to have been noticed before.
\item $\lambda$ even staircase, see Corollary \ref{escr}.  If  $\lambda=(2m,2m-2,\dots,2)$, then  $Q_\lambda(1,q,\dots,q^{n-1})$ factors. In contrast to the previously mentioned cases, we have not been able to reduce this  to previously known results.
\end{itemize}

The plan of the paper is as follows. Section \ref{ps} contains preliminaries. In Section~\ref{ks} we consider the limit $n\rightarrow \infty$ and in Section \ref{dss} the odd staircase, describing the connection to previously known results. After this preliminary material, we lay the foundation for our main results in Theorem \ref{pkt}, where we relate  $Q_\lambda(1,q,\dots,q^{n})$ to
 multivariable Christoffel--Darboux kernels, obtaining as a consequence new pfaffian and determinantal identities. These are then used  to derive the pfaffian hypergeometric identity of  Theorem \ref{th} and the  determinantal
  hypergeometric identities in Section \ref{ds}. Finally, in the Appendix we give an alternative proof of
Theorem~\ref{pkt}.

{\bf  Acknowledgement:} I  thank 
John Stembridge for several useful comments. 

\section{Preliminaries}
\label{ps}

We will use the terminology and notation of \cite{m} for partitions
and symmetric functions, and that of \cite{gr} for classical  and
basic hypergeometric functions.

\subsection{Schur $Q$-polynomials and marked shifted tableaux}
\label{tss}

When  $m\leq n$ and 
$\lambda=(\lambda_1,\dots,\lambda_m)$ is a  partition of length $m$, that is,
$$\lambda_1\geq\lambda_2\geq\dots\geq\lambda_m> 0,$$
the \emph{Schur $Q$-polynomial} $Q_\lambda$ is defined by 
\begin{equation}\label{pdef}Q_\lambda(x_1,\dots,x_n)=2^m\sum_{\sigma\in
  S_n/S_{n-m}}\sigma\Bigg(x_1^{\lambda_1}\dotsm
  x_m^{\lambda_m}\prod_{\substack{1\leq i\leq m\\1\leq i<j\leq
  n}}\frac{x_i+x_j}{x_i-x_j}\Bigg).\end{equation}
Here, $S_n$ acts by permuting the variables $x_1,\dots,x_n$ and 
$S_{n-m}$ is the  subgroup acting on $x_{m+1},\dots,x_n$. 
Note that $Q_\lambda=0$ unless $\lambda$ is strict, that is,
$$\lambda_1>\lambda_2>\dots>\lambda_m> 0.$$

For strict  $\lambda$,  $Q_\lambda$ has a combinatorial interpretation 
in terms of  \emph{marked shifted
tableaux}. We briefly recall the relevant definitions.

Let $\lambda$ be the 
diagram of a strict partition, and let $S(\lambda)$ denote the diagram obtained by shifting the
$i$-th row  $(i-1)$ steps to the right, for each $i$. A marked shifted
tableau of shape $S(\lambda)$ is a labelling of the boxes of
$S(\lambda)$ with symbols from the ordered alphabet $1'<1<2'<2<\dotsm$
such that:
 
(1) The labels increase weakly along rows and down columns.

(2) Each unmarked symbol occurs at most once in each column.

(3) Each marked symbol occurs at most once in each row.

Let $a_k$ be 
the number of boxes labelled either $k$
or $k'$, and let $x^T$ denote the monomial
$\prod_{k\geq 1}
x_k^{a_k}$. Then,  \cite[(III.8.16$'$)]{m} 
\begin{equation}\label{mcf}Q_\lambda(x)=\sum_T x^{T},\end{equation}  
where the sum is over all marked shifted tableaux of shape
$S(\lambda)$.  In \eqref{mcf} the number
of variables is  infinite; restricting the alphabet to
$1'<1<\dots<n'<n$ gives a formula for $Q_\lambda(x_1,\dots,x_n)$.

When $T$ is  a marked shifted tableau, we let $|T|=\sum_{k\geq
  1}(k-1)\,a_k$, with $a_k$ as above. Then,
\begin{equation}\label{qe}Q_\lambda(1,q,\dots,q^{n-1})=\sum_T q^{|T|} 
\end{equation}
is the generating function for $|T|$ on  marked shifted
tableaux of shape $\lambda$ and alphabet $1'<1<\dots<n'<n$.
In particular, the cardinality of this set is
$Q_\lambda(1^n)$.

\subsection{Schur $Q$-polynomials and basic hypergeometric series}

Representing  $\sigma$ by the $m$-tuple
  $(k_1,\dots,k_m)=(\sigma(1),\dots,\sigma(m))$, we may rewrite  \eqref{pdef}
   as 
\begin{align}
&Q_\lambda(x_1,\dots,x_n)\notag\\
&\quad=2^m\sum_{\substack{1\leq k_1,\dots,k_m\leq n\\k_i\text{ distinct}}}
\,\prod_{i=1}^m x_{k_i}^{\lambda_i}\prod_{1\leq i<j\leq
  m}\frac{x_{k_i}+x_{k_j}}{x_{k_i}-x_{k_j}}\prod_{i=1}^m\prod_{j\in\{1,\dots,n\}\setminus\{k_1,\dots,k_m\}}\frac{x_{k_i}+x_j}{x_{k_i}-x_j}\notag\\
&\quad=2^m\sum_{k_1,\dots,k_m=1}^n\,\prod_{i=1}^m
x_{k_i}^{\lambda_i}\prod_{1\leq i<j\leq
  m}\frac{x_{k_j}-x_{k_i}}{x_{k_j}+x_{k_i}}
\prod_{i=1}^m\prod_{j=1,\,j\neq k_i}^n\frac{x_{k_i}+x_j}{x_{k_i}-x_j}.\label{pda}
\end{align}
Consider the case when  $x_1,\dots,x_n$ are
in geometric progression. Then, \eqref{pda} is a
multiple basic hypergeometric sum. Indeed,  
replacing $n$ by $n+1$ and $k_i$ by $k_i+1$,  we obtain after simplification
\begin{multline}\label{phg}Q_\lambda(1,q,\dots,q^n)\\
=\frac{(-1;q)_{n+1}^m}
{(q;q)_n^m}
\sum_{k_1,\dots,k_m=0}^n\,\prod_{1\leq i<j\leq
  m}\frac{q^{k_j}-q^{k_i}}{q^{k_j}+q^{k_i}}\prod_{i=1}^m\frac{(-q,q^{-n};q)_{k_i}}{(q,-q^{-n};q)_{k_i}}\,q^{\lambda_ik_i},\end{multline}
 where we use the standard notation \cite{gr}
$$(a;q)_k=(1-a)(1-aq)\dotsm(1-aq^{k-1}), $$
$$(a_1,\dots,a_m;q)_k=(a_1;q)_k\dotsm(a_m;q)_k.$$

It will be convenient to view  $q^{\lambda_i}$ as free
variables. Thus, we introduce the polynomials
\begin{equation}\label{pd}P_n(x_1,\dots,x_m)=\frac{(-q;q)_n^m}{(q;q)_n^m}
\sum_{k_1,\dots,k_m=0}^n\,\prod_{1\leq i<j\leq
  m}\frac{q^{k_j}-q^{k_i}}{q^{k_j}+q^{k_i}}\prod_{i=1}^m\frac{(-q,q^{-n};q)_{k_i}}{(q,-q^{-n};q)_{k_i}}\,x_i^{k_i},\end{equation}
so that \eqref{pq} holds.
By anti-symmetry, we may equivalently write
\begin{multline}\label{qse}P_n(x_1,\dots,x_m)\\
=\frac{(-q;q)_n^m}{(q;q)_n^m}
\sum_{0\leq k_m<\dots<k_1\leq n}\,\prod_{1\leq i<j\leq
  m}\frac{q^{k_j}-q^{k_i}}{q^{k_j}+q^{k_i}}\prod_{i=1}^m\frac{(-q,q^{-n};q)_{k_i}}{(q,-q^{-n};q)_{k_i}}\,\det_{1\leq i,j\leq m}\big(x_i^{k_j}\big).\end{multline}
This can be viewed as a Schur polynomial expansion, see \eqref{pse}. One easily verifies that
\begin{equation}\label{peo}P_n(x_1,\dots,x_m,0)=(-1)^m\frac{(-q;q)_n}{(q;q)_n}\, x_1\dotsm x_m\, P_{n-1}(x_1,\dots,x_m).\end{equation}

For $m=1$, $P_n$ is a  terminating well-poised
${}_2\phi_1$ series, and can thus be expressed in terms of continuous
 $q$-ultraspherical polynomials
\cite{gr} as
\begin{equation}\label{ovp}P_n(x)=\frac{(-q;q)_n}{(q;q)_n}\,{}_2\phi_1\left[\begin{matrix}-q,q^{-n}\\-q^{-n}\end{matrix};q,x\right]=e^{in\theta}C_n(\cos\theta;-q|q),\qquad
x=-e^{2i\theta}.  
 \end{equation}
It will be convenient to introduce the monic polynomials
\begin{equation}\label{mup}c_n(x)=\frac{(q;q)_n}{(-q;q)_n}\,C_n(x/2;-q|q).
\end{equation}
Then, assuming a fixed choice of $\sqrt{-x}$,
\begin{equation}\label{pmu}P_n(x)=\frac{(-q;q)_n}{(q;q)_n}\,\left(\sqrt{-x}\right)^n c_n\left(\sqrt {-x}+({\sqrt {-x}})^{-1}\right). \end{equation}
For $0<q<1$, the polynomials $(c_k)_{k=0}^\infty$ form an orthogonal system with respect to a unique positive measure  (see \eqref{oi}), which we can normalize so that
\begin{equation}\label{cn}\| c_k\|^2      =\frac{(q;q)_k(q;q)_{k+1}}{(-q;q)_k(-q;q)_{k+1}}.\end{equation}

\subsection{Orthogonal polynomials}

Our main interest is in the case $q=1$. Unfortunately,  
the right-hand side of \eqref{phg} is then singular. However,
 when
 $m=1$ one can use  transformation formulas for basic hypergeometric
 series to derive a plethora of other expressions, many
of which remain regular when $q\rightarrow 1$. 

The transformations that we will use correspond to different identifications of $q$-ultraspherical polynomials with continuous $q$-Jacobi polynomials, which in turn form a sub-class of the 
  Askey--Wilson polynomials. Recall that the latter polynomials are defined by
$$p_n(\cos\theta ;a,b,c,d|q)=\frac{(ab,ac,ad;q)_n}{a^n}\,{}_4\phi_3\left[\begin{matrix}q^{-n},abcdq^{n-1},ae^{i\theta},ae^{-i\theta}\\ab,ac,ad\end{matrix};q,q\right]; $$
they are symmetric in the  parameters $a,b,c,d$ and satisfy
\begin{equation}\label{awp}p_n(-x;a,b,c,d|q)=(-1)^np_n(x;-a,-b,-c,-d|q). \end{equation}
 By  \cite[Eq.\ (7.5.33--36)]{gr},
\begin{subequations}\label{cpid}
\begin{align}C_n(x;a|q)&=\frac{(a;q)_n}{(q,a^2q^n;q)_n}\,p_n(x;\sqrt a,-\sqrt a, \sqrt{aq},-\sqrt{aq}|q)\\
 &=\frac{(a; q^{1/2})_n}{(q,aq^{1/2};q)_n}\,p_n(x;\sqrt a,-\sqrt a, q^{1/4},-q^{1/4}|q^{1/2}),
\end{align}
\begin{align}C_{2n}(x;a|q)&=\frac{(a;q)_n}{(q;q)_{2n}}\,p_n(2x^2-1;a,-1,\sqrt q,-\sqrt q|q)\\
&=\frac{(a^2;q^2)_n}{(q,-a;q)_{2n}}\,p_n(2x^2-1;a,a q,-1,-q|q^2),\end{align}
\begin{align}C_{2n+1}(x;a|q)&=
\frac{(a;q)_{n+1}}{(q;q)_{2n+1}}\,2x\,p_n(2x^2-1;a,-q,\sqrt q,-\sqrt q|q)
\\
&=\frac{(a^2;q^2)_{n+1}}{(q,-a;q)_{2n+1}}\,2x\,p_n(2x^2-1;a,a q,-q,-q^2|q^2).
\end{align}
\end{subequations}
It is  easy to prove this by verifying that the orthogonality measures for the various  polynomials agree.

In view of \eqref{ovp}, we are particularly interested in the case  $a=-q$. Using also  \eqref{awp}, we may then write \eqref{cpid} as 
\begin{subequations}\label{cbid}
\begin{align}\label{cid}P_n(-e^{2i\theta})&=\frac{(-q;q)_n}{(q,q^{n+2};q)_n}\,e^{in\theta} \,p_n(\cos\theta;i q^{1/2},-i q^{1/2}, iq,-iq|q)\\
\label{cidb} &=\frac{(-q; q^{1/2})_n}{(q,-q^{3/2};q)_n}\,e^{in\theta}\,p_n(\cos\theta;i q^{1/2},-i q^{1/2}, q^{1/4},-q^{1/4}|q^{1/2}),
\end{align}
\end{subequations}
\begin{subequations}\label{chi}
\begin{align}P_{2n}(x)&=\frac{(-q;q)_n}{(q;q)_{2n}}\,x^n\,p_n({\textstyle\frac12}(x+x^{-1});1,q,\sqrt q,-\sqrt q|q)\\
&=\frac{(q^2;q^2)_n}{(q;q)^2_{2n}}\,x^n\,p_n({\textstyle\frac12}(x+x^{-1})
;1,q,q,q^2|q^2),\end{align}
\begin{align}P_{2n+1}(x)&=
\frac{(-q;q)_{n+1}}{(q;q)_{2n+1}}\,x^n(1-x)\,p_n({\textstyle\frac12}(x+x^{-1});q,q,\sqrt q,-\sqrt q|q)
\\
&=\frac{(q^2;q^2)_{n+1}}{(q;q)^2_{2n+1}}\,x^n(1-x)\,p_n({\textstyle\frac12}(x+x^{-1});q,q,q^2,q^2|q^2).
\end{align}
\end{subequations}

These identities yield
 many useful  expressions
 for $P_n(x)$. For instance, \eqref{cid} implies
\begin{equation}\label{43}P_n(x)=\frac{1-q^{n+1}}{1-q}\left(\sqrt{x/q}\right)^{n}\,{}_4\phi_3\left[\begin{matrix}q^{-n},q^{n+2},\sqrt{q/x},-\sqrt{qx}\\q,q^{3/2},-q^{3/2}\end{matrix};q,q\right].\end{equation}
In contrast to \eqref{ovp}, this expression is regular in the limit
\begin{equation}\label{3f2}Q_{(\lambda)}(1^{n+1})=\lim_{q\rightarrow 1} 2P_n(q^{\lambda})
=(2n+2)\,{}_3F_2\left[\begin{matrix}-n,{n+2},(1-\lambda)/2\\1,3/2\end{matrix}\,;1\right].
 \end{equation}
This will be generalized to  general $m$ in Theorem~\ref{th} and Corollary \ref{qc}. 

Similarly, \eqref{cidb} gives
$$P_n(x)=\frac{1-q^{n+1}}{1-q}\left(\sqrt{x/q}\right)^{n}
\,{}_4\phi_3\left[\begin{matrix}q^{-n/2},-q^{(n+2)/2},\sqrt{q/x},-\sqrt{qx}\\q,iq^{3/4},-iq^{3/4}\end{matrix};q^{1/2},q^{1/2}\right],$$
which implies  the simple  formula
\begin{equation}\label{2f1}Q_{(\lambda)}(1^{n+1})=\lim_{q\rightarrow 1} 2P_n(q^{\lambda})
=(2n+2)\,{}_2F_1\left[\begin{matrix}-n,1-\lambda\\2\end{matrix}\,;2\right].
 \end{equation}

In Section \ref{ds}, we will need the following eight  expressions, which are all consequences of  \eqref{chi}:  
\begin{subequations}\label{chia}
\begin{align}\label{chiaa}
 P_{2n}(x)
&=\frac{1-q^{2n+1}}{1-q}\,q^{-n}x^n\,{}_4\phi_3\left[\begin{matrix}q^{-n},-q^{n+1},qx,{q/x}\\q,q^{3/2},-q^{3/2}\end{matrix};q,q\right]\\
&=\frac{(q^{3/2},-q;q)_n}{(q,-q^{1/2};q)_n}\,q^{-\frac n2}x^n\,{}_4\phi_3\left[\begin{matrix}q^{-n},-q^{n+1},q^{1/2}x,{q^{1/2}/x}\\q^{1/2},q^{3/2},-q\end{matrix};q,q\right]\\
&=\frac{1-q^{2n+1}}{1-q}\,q^{-n}x^n\,{}_4\phi_3\left[\begin{matrix}q^{-2n},q^{2n+2},{qx},q/x\\q,q^{2},q^{3}\end{matrix};q^2,q^2\right]\\
&=\left(\frac{1-q^{2n+1}}{1-q}\right)^2q^{-2n}x^n\,{}_4\phi_3\left[\begin{matrix}q^{-2n},q^{2n+2},q^2x,{q^2/x}\\q^2,q^{3},q^{3}\end{matrix};q^2,q^2\right],
\end{align}
\begin{align}
&\quad P_{2n+1}(x)
=\frac{(1-q^{2n+2})(1-x)}{(1-q)^2}\,q^{-n}x^n\,{}_4\phi_3\left[\begin{matrix}q^{-n},-q^{n+2},qx,{q/x}\\q^2,q^{3/2},-q^{3/2}\end{matrix};q,q\right]\\
&=\frac{(q^{3/2};q)_n(-q;q)_{n+1}(1-x)}{(1-q)(q,-q^{3/2};q)_n}\,q^{-\frac n2}x^n\,{}_4\phi_3\left[\begin{matrix}q^{-n},-q^{n+2},q^{1/2}x,{q^{1/2}/x}\\q^{3/2},q^{3/2},-q\end{matrix};q,q\right]\\
&=\frac{(1-q^{2n+2})(1-x)}{(1-q)^2}\,q^{-n}x^n\,{}_4\phi_3\left[\begin{matrix}q^{-2n},q^{2n+4},{qx},q/x\\q^2,q^{3},q^{3}\end{matrix};q^2,q^2\right]\\
&=\frac{(1-q^{2n+2})^2(1-x)}{(1-q)^2(1-q^2)}\,q^{-2n}x^n\,{}_4\phi_3\left[\begin{matrix}q^{-2n},q^{2n+4},q^2x,{q^2/x}\\q^3,q^{3},q^{4}\end{matrix};q^2,q^2\right].
\end{align}
\end{subequations}

In the limit $q\rightarrow 1$, the Askey--Wilson polynomials in \eqref{cbid} and \eqref{chi}  
degenerate to four different systems in the Askey Scheme. Let us elaborate on this point, introducing the notation
$$p_k(x)=
\frac{i^kk!}{2^k}\lim_{q\rightarrow 1}P_k(q^{-ix}), $$
where the constant is chosen so as to make $p_k$ monic.
By \eqref{3f2} and \eqref{2f1},
$$p_k(x)=\frac{i^k(k+1)!}{2^k}\,{}_2F_1\left[\begin{matrix}-k,1+ix\\2\end{matrix}\,;2\right]=\frac{i^k(k+1)!}{2^k}\,{}_3F_2\left[\begin{matrix}-n,{n+2},(1+ix)/2\\1,3/2\end{matrix}\,;1\right].$$
It follows that $p_k$
may be identified with Meixner-Pollaczek polynomials and with continuous Hahn polynomials. Explicitly, in the notation of  \cite{ks},
$$p_k(x)=\frac{k!}{2^k}\,P^{(1)}_k(x;\pi/2)=\frac{2^k k!(k+1)!}{(2k+1)!}\,P_k(x/2;1/2,1,1/2,1). $$
In particular, $(p_k(x))_{k=0}^\infty$ are the monic orthogonal polynomials corresponding to the measure
$$\int_{-\infty}^\infty \frac{x}{\sinh(\pi x)}\,f(x)\,dx, $$ 
with norms
$$ \|p_k\|^2=\frac{k!(k+1)!}{2^{2k+1}}.$$

Let
\begin{equation}\label{cdh}p_{2k}(x)=p_k^{(0)}(x^2),\qquad p_{2k+1}(x)=xp_k^{(1)}(x^2),\end{equation}
so that $(p_k^{(\ep)}(x))_{k=0}^\infty$ are the monic orthogonal polynomials corresponding to 
$$\int_{0}^\infty \frac{x^\ep f(x)}{\sinh(\pi \sqrt x)}\,dx,\qquad \ep=0,\,1, $$
with norms
\begin{equation}\label{cdhn}
\|p_k^{(0)}\|^2=\frac{(2k)!(2k+1)!}{2^{4k+1}},\qquad \|p_k^{(1)}\|^2=\frac{(2k+1)!(2k+2)!}{2^{4k+3}}.
\end{equation} 
Then,  $p_k^{(\ep)}$  are continuous dual Hahn and Wilson polynomials; explicitly,
$$p_k^{(0)}(x)=(-1)^k S_k(x;0,1/2,1)=(-1)^k\frac{4^kk!}{(2k)!}\,W_k(x/4;0,1/2,1/2,1), $$
$$p_k^{(1)}(x)=(-1)^k S_k(x;1/2,1,1)=(-1)^k\frac{4^k(k+1)!}{(2k+1)!}\,W_k(x/4;1/2,1/2,1,1). $$
This may be seen  by letting $q\rightarrow 1$ in \eqref{chia}, or by checking that the orthogonality measures agree. 

\begin{remark}
It may seem strange that the variables $\lambda_i$, which are originally integers, have to be taken as imaginary at the support of the orthogonality measure of $p_k$. Actually, if we let
\begin{equation}\label{fk}f_k(x)=\frac{p_k(ix)}{i^k}=
\frac{(k+1)!}{2^k}\,
{}_2F_1\left[\begin{matrix}-k,1-x\\2\end{matrix}\,;2\right], \end{equation}
then   $f_k$ are ``almost'' orthogonal  on the positive integers, in view of the Abel sum
$$\lim_{t\rightarrow 1} \sum_{k=1}^\infty (-1)^{k+1}t^kk\,f_m(k)f_n(k)
=\frac{(-1)^n(n+1)!n!}{4^{n+1}}\,\delta_{mn}.
 $$
This type of orthogonality plays an important role in \cite{rss}.
It is a limit case of the orthogonality of Meixner polynomials $M_n(x;\beta,c)$ \cite{ks} with $\beta=2$, $c\rightarrow -1$. 
When $q\neq 1$, a similar formula can be obtained  as a limit case of the discrete orthogonality for $q$-ultraspherical polynomials considered in \cite{ak}.
\end{remark}

\subsection{Schur $Q$-polynomials and pfaffians}

The \emph{pfaffian} of a skew-symmetric even-dimensional matrix
is defined by
$$\pfaff_{1\leq i,j\leq 2m}(a_{ij})
=\sum_{\sigma\in S_{2m}/G}\sgn(\sigma)\prod_{i=1}^m
a_{\sigma(2i-1),\sigma(2i)},$$
where $G$ is the subgroup of order $2^m m!$ consisting of
 permutations preserving the set of  pairs 
$\{\{1,2\},\{3,4\},\dots,\{2m-1,2m\}\}$. Pfaffians and Schur
$Q$-polynomials are intimately related. Indeed, 
Schur's original definition of the latter \cite{s} is based on the
identity
\begin{equation}\label{sqp}
Q_{(\lambda_1,\dots,\lambda_{2m})}=\pfaff_{1\leq i,j\leq
  2m}\left(Q_{(\lambda_i,\lambda_j)}\right). \end{equation}
In the same work, Schur obtained the pfaffian evaluation
\begin{equation}\label{spa}\pfaff_{1\leq i,j\leq
  2m}\left(\frac{x_j-x_i}{x_j+x_i}\right)=\prod_{1\leq i<j\leq
  {2m}}\frac{x_j-x_i}{x_j+x_i},\end{equation}
which we recall here together with its companion 
\cite[Proposition 2.3]{ste} 
\begin{equation}
\label{spb}\pfaff_{1\leq i,j\leq
  2m}\left(\frac{x_j-x_i}{1-tx_ix_j}\right)=t^{m(m-1)}\prod_{1\leq i<j\leq
  {2m}}\frac{x_j-x_i}{1-tx_ix_j}.
\end{equation}

We  need an elementary property of pfaffians, which  is
 stated below as  
 Corollary \ref{lms}. Incidentally, it can be used to prove 
 \eqref{sqp}, see Remark \ref{sqr}.
 Before formulating the result, we find it instructive to give a generalization. 

\begin{lemma}\label{dbl}
Let $(X_i,\mu_i)_{1\leq i\leq 2m}$ be a collection of measure spaces, and 
$(b_{ij})_{1\leq i,j\leq 2m}$ a collection of integrable functions, $b_{ij}:\,X_i\times X_j\rightarrow \mathbb R$, such that
$b_{ij}(x,y)=-b_{ji}(y,x)$. Then,
\begin{multline*}\pfaff_{1\leq i,j\leq 2m}\left(\iint b_{ij}(x,y)\,d\mu_i(x)d\mu_j(y)\right) \\
=\idotsint\pfaff_{1\leq i,j\leq 2m}\left( b_{ij}(x_i,x_j)\right)d\mu_1(x_1)\dotsm d\mu_m(x_{2m}).\end{multline*}
\end{lemma}

\begin{proof}
The left-hand side  equals
$$\frac{1}{2^m m!}\sum_{\sigma\in S_{2m}}\sgn(\sigma)\prod_{i=1}^m
\iint b_{\sigma(2i-1),\sigma(2i)}(x_i,y_i)\,d\mu_{\sigma(2i-1)}(x_i)d\mu_{\sigma(2i)}(y_i).$$
Introducing new integration variables by $x_i\mapsto x_{\sigma(2i-1)}$, 
$y_i\mapsto x_{\sigma(2i)}$ and interchanging the finite sum and the integral, we obtain
 the desired result.
\end{proof}

Consider the special case of Lemma \ref{dbl} when each $X_i$ is equal to the same 
finite  space $X$. Writing
$X=\{1,\dots,n\}$, $\mu_i(j)=A_{ij}$, $b_{ij}(k,l)=B^{ij}_{kl}$, it takes the following form.

\begin{corollary}\label{lms}
Let $(A_{ij})_{1\leq i\leq 2m,\,1\leq j\leq n}$ and $(B_{ij}^{kl})_{1\leq i,j\leq
  n,\,1\leq k,l\leq 2m}$ be matrices such that $B^{kl}_{ij}=-B^{lk}_{ji}$. Then, 
$$\pfaff_{1\leq i,j\leq
    2m}\left(\sum_{x,y=1}^nA_{ix}A_{jy} B^{ij}_{xy}\right)
=\sum_{k_1,\dots,k_{2m}=1}^n\,\prod_{i=1}^{2m}
  A_{ik_i}\pfaff_{1\leq i,j\leq 2m}(B_{k_ik_j}^{ij}).$$
 \end{corollary}

\begin{remark}
We will only need Corollary \ref{lms} in the case when 
 $B_{ij}=B_{ij}^{kl}$ is independent of $k$ and $l$. It can then
  equivalently be written
$$\pfaff(ABA^{ t})=\sum_{1\leq k_{2m}<\dots<k_{1}\leq n}
\det_{1\leq i,j\leq 2m}(A_{i,k_j})\pfaff_{1\leq i,j\leq 2m}(B_{k_ik_j}),
$$
which is a version of the ``minor summation formula'' of Ishikawa and Wakayama
 \cite{iw1}. The corresponding special case of Lemma \ref{dbl} 
 is due to de Bruijn \cite{b}. 
\end{remark}

\begin{remark}\label{sqr}
To prove \eqref{sqp}, one may apply
 Corollary \ref{lms} in the case when
$$A_{ik}=x_k^{\lambda_i}\prod_{j=1,\,j\neq k}^n\frac{x_k+x_j}{x_k-x_j},\qquad B_{ij}^{kl}=B_{ij}=\frac{x_j-x_i}{x_j+x_i}. $$
Using  \eqref{spa}
to compute $\pfaff(B_{k_ik_j})$
and comparing with \eqref{pda} gives indeed \eqref{sqp}. 
\end{remark}

Similarly, choosing
$$A_{ik}=x_i^{k-1}\frac{(-q,q^{-n};q)_{k-1}}{(q,-q^{-n};q)_{k-1}},\qquad B_{ij}^{kl}=B_{ij}=\frac{q^j-q^i}{q^j+q^i}. $$
it follows from \eqref{pd} that, for $m$ even,
\begin{equation}\label{pfi}
P_n(x_1,\dots,x_{m})=\pfaff_{1\leq i,j\leq m}\left(P_n(x_i,x_j)\right).\end{equation}
Together with \eqref{peo}, this reduces the study of the polynomials 
$P_n(x_1,\dots,x_{m})$ to the case $m=2$.

\subsection{Multivariable Christoffel--Darboux kernels}

Let 
\begin{equation}\label{mf}f\mapsto \int f(x)\,d\mu(x)\end{equation}
be a positive moment functional and
$(p_k(x))_{k=0}^\infty$ the corresponding family of monic orthogonal polynomials.
We may then introduce the Christoffel--Darboux kernel
 \begin{equation}\label{cd}K^n(x,y)=\sum_{k=0}^{n-1}\frac{p_k(x)p_k(y)}{\| p_k\|^2}
=\frac 1{\| p_{n-1}\|^2}\frac{p_n(x)p_{n-1}(y)-p_{n-1}(x)p_n(y)}{x-y}.\end{equation}

In \cite{r0}, we studied  more general multivariable kernels 
$$K_m^n(x)=K_m^n(x_1,\dots,x_{2m}).$$
For  present purposes,
  the crucial fact is that the following three expressions define the same kernel.  As is explained in \cite{r0}, this can be deduced from results of  Ishikawa and Wakayama \cite{iw2}, Lascoux \cite{l1,l2}, Okada \cite{o} and Strahov and Fyodorov \cite{sf}. A self-contained  proof is given in \cite{r0}.

\begin{lemma}\label{cdl}
If
\begin{equation}\label{kd}K_{m}^n(x)=\frac{1}{\prod_{i=1}^m\| p_{n-i}\|^2}
\frac{\det_{1\leq i,j\leq 2m}(p_{n-m+j-1}(x_i))}{\prod_{1\leq i<j\leq 2m}(x_j-x_i)}, 
\end{equation}
then, for any choice of square roots  $\sqrt {x_i}$,
$$K_m^n(x)=\frac {1}{\prod_{1\leq i<j\leq 2m}(\sqrt {x_j}-\sqrt {x_i})}\pfaff_{1\leq i,j\leq 2m}\left((\sqrt {x_j}-\sqrt {x_i})K(x_i,x_j)\right) $$
and, for any choice of $\xi_i$ such that  $\xi_i+\xi_i^{-1}=x_i+2$,
$$K_m^n(x)=\frac {\prod_{i=1}^{2m} \xi_i^{m-1}}{\prod_{1\leq i<j\leq 2m}( {\xi_j}- {\xi_i})}\pfaff_{1\leq i,j\leq 2m}\left(( {\xi_j}- {\xi_i})K(x_i,x_j)\right). $$
\end{lemma}

We note in passing the less symmetric identities
\begin{multline}\label{mk}K_{m}^n(x_1,\dots,x_m,y_1,\dots,y_m)=\frac{\det_{1\leq i,j\leq m}(K^n(x_i,y_j))}
{\prod_{1\leq i<j\leq m}(x_j-x_i)(y_j-y_i)}\\
\shoveleft{=\frac{1}{\prod_{1\leq i<j\leq m}(x_j-x_i)(y_j-y_i)}}\\
\times\sum_{0\leq k_m<\dots<k_1\leq n-1}\,
\prod_{i=1}^m\frac{1}{\|p_{k_i}\|^2}\,\det_{1\leq i,j\leq m}(p_{k_i}(x_j))\det_{1\leq i,j\leq m}(p_{k_i}(y_j)). \end{multline}
Thus, we have two determinant and two pfaffian formulas for the 
 multiple sum in \eqref{mk}; these give four different multivariable extensions of the Christoffel--Darboux formula \eqref{cd}. 

We will apply Lemma \ref{cdl} in  slightly modified form. To this end, assume that all odd moments of the functional \eqref{mf} vanish.  Let  $(q_k)_{k=0}^\infty$ and $(r_k)_{k=0}^\infty$ be the monic orthogonal polynomials associated with the moment functionals
$$f\mapsto\int f(-x^2)\,d\mu(x),\qquad f\mapsto\int f(-x^2)x^2\,d\mu(x), $$
respectively; the minus signs will  be convenient. Then,
$$p_{2k}(x)=(-1)^kq_k(-x^2),\qquad p_{2k+1}(x)=(-1)^kx\,r_k(-x^2). $$
Later, we will choose  $p_k=c_k$ as the continuous
 $q$-ultraspherical polynomials in
\eqref{mup}. Then,
one may use \eqref{chi} to identify $q_k$ and $r_k$ with continuous $q$-Jacobi polynomials.

Let 
 \begin{equation}\label{ecd}\tilde K^n(x,y)= \sum_{\substack{0\leq k\leq n-1\\k\equiv n-1\,\operatorname{mod}\, 2}}\frac{p_k(x)p_k(y)}{\| p_k\|^2}.\end{equation}
Then, 
applying \eqref{cd} to the systems $(q_k)_{k=0}^\infty$ and $(r_k)_{k=0}^\infty$ gives 
\begin{equation}\label{ecd2}\tilde K^n(x,y)= \frac 1{\| p_{n-1}\|^2}\frac{p_{n+1}(x)p_{n-1}(y)-p_{n+1}(x)p_{n-1}(y)}{x^2-y^2}.\end{equation}
We note in passing  that \eqref{ecd} implies
$$\tilde K^n(x,y)= K^n(x,y)-\tilde K^{n-1}(x,y), $$
which, upon iteration, yields
\begin{equation}\label{ktk}\tilde K^n(x,y)=\sum_{j=0}^{n-1}(-1)^{n+j+1} K^{j+1}(x,y). \end{equation}

More generally, let us temporarily write $Q_m^n$ and $R_m^n$ for the kernels obtained by replacing $(p_k)_{k=0}^\infty$ in  $K_m^n$ with $(q_k)_{k=0}^\infty$ and $(r_k)_{k=0}^\infty$, respectively. Moreover, let
$$\tilde K_m^{2n-1}(x_1,\dots,x_{2m})=Q_m^n(-x_1^2,\dots,-x_{2m}^2),
 $$
$$\tilde K_m^{2n}(x_1,\dots,x_{2m})=x_1\dots x_{2m}R_m^n(-x_1^2,\dots,-x_{2m}^2).
 $$
Then, applying Lemma \ref{cdl} to the kernels $Q$ and $R$, the result may be expressed in unified form as follows.

\begin{corollary}\label{cdc}
Under the assumption of vanishing odd moments, the following identities hold:
\begin{eqnarray}\notag\tilde K_m^n(x)&=&\frac{1}{\prod_{i=1}^m\| p_{n+1-2i}\|^2}
\frac{\det_{1\leq i,j\leq 2m}(p_{n-2m+2j-1}(x_i))}{\prod_{1\leq i<j\leq m}(x_j^2-x_i^2)}\\
\label{kpa}&=&\frac {1}{\prod_{1\leq i<j\leq 2m}( {w_j}- {w_i})}\pfaff_{1\leq i,j\leq 2m}\left(( {w_j}- {w_i})\tilde K^n(x_i,x_j)\right)\\
\label{kpb}&=&\frac {\prod_{i=1}^{2m} \xi_i^{m-1}}{\prod_{1\leq i<j\leq 2m}( {\xi_j}- {\xi_i})}
\pfaff_{1\leq i,j\leq 2m}\left(( {\xi_j}- {\xi_i})\tilde K^n(x_i,x_j)\right),
\end{eqnarray}
where $w_i^2=-x_i^2$ and $\xi_i+\xi_i^{-1}=2-x_i^2$,
 $1\leq i\leq 2m$.
\end{corollary}

The right-hand side of \eqref{kd} is an analogue of a rectangular Schur polynomial; in fact, it can be identified with a Schur function over an abstract alphabet \cite{l1}. These polynomials make sense also for odd-dimensional matrices. Let
\begin{equation}\label{rsd}\mathbb P_{n^m}(x_1,\dots,x_m)=\frac{\det_{1\leq i,j\leq m}(p_{n+j-1}(x_i))}{\det_{1\leq i,j\leq m}(p_{j-1}(x_i))}=\frac{\det_{1\leq i,j\leq m}(p_{n+j-1}(x_i))}{\prod_{1\leq i<j\leq m}(x_j-x_i)}, \end{equation}
so that
$$K_m^n(x)=\frac 1{\prod_{i=1}^m\|p_{n-i}\|^2}\,\mathbb P_{(n-m)^{2m}}(x). $$
Under the assumption of vanishing odd moments, we will also write
\begin{equation}\label{mrs} \tilde{\mathbb P}_{n^m}(x_1,\dots,x_m)=\frac{\det_{1\leq i,j\leq m}(p_{n+2j-2}(x_i))}{\prod_{1\leq i<j\leq m}(x_j^2-x_i^2)}, \end{equation}
so that
$$\tilde K_m^n(x)=\frac 1{\prod_{i=1}^m\|p_{n+1-2i}\|^2}\,\tilde{\mathbb P}_{(n+1-2m)^{2m}}(x). $$
Then, if $\mathbb Q$ and $\mathbb R$ denote the kernels obtained by replacing $p_k$ in \eqref{rsd} with $q_k$ and $r_k$, respectively, we have
$$\tilde{\mathbb  P}_{(2n)^m}(x_1,\dots,x_{m})=(-1)^{mn}\,\mathbb Q_{n^m}(-x_1^2,\dots,-x_{m}^2),
 $$
$$\tilde{\mathbb  P}_{(2n+1)^m}(x_1,\dots,x_{m})=(-1)^{mn}x_1\dotsm x_m\,\mathbb R_{n^m}(-x_1^2,\dots,-x_{m}^2).
 $$

It will be useful to note that
\begin{equation}\label{sfl}\lim_{t\rightarrow\infty}t^{-n}\tilde{\mathbb  P}_{n^{m+1}}(x_1,\dots,x_{m},t)
=\tilde{\mathbb  P}_{n^m}(x_1,\dots,x_{m});
 \end{equation}
this is easily derived from \eqref{mrs}.

Finally, we mention the integral formula
\begin{multline*}\mathbb P_{n^m}(x_1,\dots,x_m)\\
=\frac{1}{n!\prod_{i=1}^n\|p_{i-1}\|^2}
\int\prod_{\substack{1\leq j\leq m\\1\leq k\leq n}}(x_j-y_k)\prod_{1\leq i<j\leq n}
(y_j-y_i)^2\,d\mu(y_1)\dotsm d\mu(y_{n}),
 \end{multline*}
which is at least implicitly due to Christoffel, see also \cite{bh,r0}.
Applying this to the polynomials $q_k$ and $r_k$ gives
\begin{subequations}\label{kif}
\begin{multline}\tilde{\mathbb P}_{(2n)^m}(x_1,\dots,x_m) =\frac{1}{n!\prod_{i=1}^n\|p_{2i-2}\|^2}\\
\times\int\prod_{\substack{1\leq j\leq m\\1\leq k\leq n}}(x_j^2-y_k^2)\prod_{1\leq i<j\leq n}
(y_j^2-y_i^2)^2\,d\mu(y_1)\dotsm d\mu(y_{n}),
 \end{multline}
\begin{multline}\tilde{\mathbb P}_{(2n+1)^m}(x_1,\dots,x_m)=\frac{x_1\dotsm x_{m}}{n!\prod_{i=1}^n\|p_{2i-1}\|^2}\\
\times\int\prod_{\substack{1\leq j\leq m\\1\leq k\leq n}}(x_j^2-y_k^2)\prod_{1\leq i<j\leq n}
(y_j^2-y_i^2)^2\prod_{k=1}^{n}y_k^2\,d\mu(y_1)\dotsm d\mu(y_{n}).
 \end{multline}
\end{subequations}
These formulas will be used in the Appendix.

\section{Relation to Kawanaka's identity}
\label{ks}

Kawanaka \cite[Theorem 1.1]{k} obtained the Schur polynomial identity
\begin{equation}\label{kbf}\sum_\mu
\,\prod_{\alpha\in\mu}\frac{1+q^{h(\alpha)}}{1-q^{h(\alpha)}}\,
s_\mu(x_1,\dots,x_m)=\prod_{i=1}^m\frac{(-x_iq;q)_\infty}{(x_i;q)_\infty}\prod_{1\leq
i<j\leq m}\frac{1}{1-x_ix_j}.\end{equation}
Here,  the sum is over all
partitions $\mu_1\geq\mu_2\geq\dots\geq\mu_m\geq 0$ and
$h(\alpha)$ denotes the \emph{hook-length} at the box $\alpha$ of the 
diagram of $\mu$, see \cite{m}. We will recover a proof of \eqref{kbf} below, see
Remark \ref{kr}.

Although it is not mentioned by Kawanaka,
 \eqref{kbf} can be viewed as evaluating a Schur $Q$-function at an
 \emph{infinite} geometric progression. To see this we
 observe that, 
writing   $k_i=\mu_i+m-i$, the determinant in \eqref{qse} equals
$$\prod_{1\leq i<j\leq m}(x_i-x_j)\,s_\mu(x_1,\dots,x_m). $$
Moreover, 
it is easy to check that (cf.\ \cite[Examples I.1.1 and I.1.3]{m})
$$\prod_{i=1}^m\frac{(-q;q)_{k_i}}{(q;q)_{k_i}}\prod_{1\leq i<j\leq
  m}\frac{q^{k_j}-q^{k_i}}{q^{k_j}+q^{k_i}}=\prod_{\alpha\in\mu}\frac{1+q^{h(\alpha)}}{1-q^{h(\alpha)}},
  $$
$$\prod_{i=1}^m\frac{(q^{-n};q)_{k_i}}{(-q^{-n};q)_{k_i}}
=(-1)^{\binom m2}\prod_{i=1}^m\frac{(q^{n+2-i};q)_{i-1}}{(-q^{n+2-i};q)_{i-1}}
\prod_{\alpha\in\mu}\frac{1-q^{c(\alpha)+m-n-1}}{1+q^{c(\alpha)+m-n-1}},
 $$
where $c$ denotes  \emph{content}.
Plugging all this into \eqref{qse} gives
\begin{multline}\label{pse}
P_n(x_1,\dots,x_m)=\prod_{i=1}^m
\frac{(-q;q)_{n+1-i}}{(q;q)_{n+1-i}}
\prod_{1\leq i<j\leq
  m}(x_j-x_i)\\
\times\sum_{0\leq\mu_m\leq\dots\leq
  \mu_1\leq
  n+1-m}\,\prod_{\alpha\in\mu}\frac{1+q^{h(\alpha)}}{1-q^{h(\alpha)}}\frac{1-q^{c(\alpha)+m-n-1}}{1+q^{c(\alpha)+m-n-1}}\,s_\mu(x). 
\end{multline}

Let us now take the  limit
 $n\rightarrow\infty$, which can be justified
 analytically  if $|q|<1$ and $|x_i|<1$ for all $i$. Note that
$$\lim_{n\rightarrow\infty}\prod_{\alpha\in\mu}\frac{1-q^{c(\alpha)+m-n-1}}{1+q^{c(\alpha)+m-n-1}}=(-1)^{\sum_i\mu_i},
 $$
which may be absorbed into $s_\mu$ by homogeneity. Thus, we obtain
 \begin{equation}\label{k2}
\lim_{n\rightarrow\infty}P_n(x_1,\dots,x_m)=\frac{(-q;q)_\infty^m}{(q;q)_\infty^m}
\prod_{1\leq i<j\leq
  m}(x_j-x_i)\sum_{\mu}\,\prod_{\alpha\in\mu}\frac{1+q^{h(\alpha)}}{1-q^{h(\alpha)}}\,s_\mu(-x), 
\end{equation}
where the sum  is computed by \eqref{kbf}.
Although the corresponding Schur $Q$-function
identity is equivalent to \eqref{kbf}, it 
seems not to have appeared explicitly in the literature.

\begin{proposition}\label{wsp}
For $\lambda$ a partition of length $m$, 
$$
Q_{\lambda}(1,q,q^2,\dotsc)=\prod_{i=1}^m\frac{(-1;q)_{\lambda_i}}{(q;q)_{\lambda_i}}
\prod_{1\leq i<j\leq
  m}\frac{q^{\lambda_j}-q^{\lambda_i}}{1-q^{\lambda_i+\lambda_j}},
$$
where the left-hand side is interpreted as
$$\lim_{n\rightarrow\infty}Q_{\lambda}(1,q,\dots,q^n),\qquad |q|<1. $$
\end{proposition}

In terms of tableaux, the result may be written
 as follows.

\begin{corollary}\label{kpp}
For $\lambda$ a strict partition of length $m$, 
$$
\sum_T q^{|T|}=\prod_{i=1}^m\frac{(-1;q)_{\lambda_i}}{(q;q)_{\lambda_i}
}
\prod_{1\leq i<j\leq
  m}\frac{q^{\lambda_j}-q^{\lambda_i}}{1-q^{\lambda_i+\lambda_j}},\qquad |q|<1,
$$
where the sum is over all marked shifted tableaux of shape
$S(\lambda)$. 
\end{corollary}

This should be compared with the identity
\begin{equation}\label{gi}\sum_T q^{|T|}=\prod_{i=1}^m\frac{1}{(q;q)_{\lambda_i}
}
\prod_{1\leq i<j\leq
  m}\frac{q^{\lambda_j}-q^{\lambda_i}}{1-q^{\lambda_i+\lambda_j}}, 
\end{equation}
where the sum is over column-strict shifted tableaux of shape
$S(\lambda)$, that is,
  only unmarked symbols appear, so that  conditions (1)--(3) in Section
 \ref{tss} may
 be summarized as 

(4) The labels increase weakly along rows and strictly down columns.

As was pointed out by Stembridge \cite{ste}, \eqref{gi} is equivalent
to a result conjectured by Stanley \cite{sta} and first proved by
Gansner \cite{g}. Comparison with
Corollary~\ref{kpp} suggests the following  problem.

\begin{problem}\label{p}
Prove directly (that is, without summing the series)
that
$$\sum_{T\text{ {marked}}}q^{|T|}=\prod_{i=1}^m(-1;q)_{\lambda_i}\sum_{T\text{ {column-strict}}} q^{|T|}, $$
where $\lambda$ is a strict partition of length $m$ and the sums are over shifted tableaux of shape $S(\lambda)$.
For instance, writing
 $M_k$ and $C_k$ for the set of marked, respectively column-strict,
shifted tableaux $T$ of shape $S(\lambda)$ such that $|T|=k$, construct a  map
$\phi:\,M_k\rightarrow\bigcup_{j\leq k} C_j$ such that
$$ \sum_k|M_{k+j}\cap\phi^{-1}(x)|\,
q^k=\prod_{i=1}^m(-1;q)_{\lambda_i},\qquad x\in C_j.
$$
This would lead to a new proof of the Stanley--Gansner identity as a
consequence of Kawanaka's identity (and vice versa). 
\end{problem}

\section{Relation to a discrete Selberg integral}
\label{dss}

Krattenthaler \cite[Theorem 6]{kr} gave the following
 multiple extension of
the $q$-Chu--Vandermonde summation, 
which was applied to  enumeration problems for perfect
matchings:
\begin{multline}\label{krs}
\sum_{0\leq k_m<\dots<k_1\leq n}\,\prod_{1\leq i<j\leq
  m}(q^{k_j}-q^{k_i})^2
\prod_{i=1}^m\frac{(x;q)_{k_i}(y;q)_{n-k_i}}{(q;q)_{k_i}(q;q)_{n-k_i}}\,y^{k_i}\\
=q^{2\binom m3}y^{\binom m2}\prod_{i=1}^m\frac{(x,y,q;q)_{i-1}(xyq^{i+m-2};q)_{n+1-m}}{(q;q)_{n+i-m}}.
\end{multline}
An equivalent identity was  
previously obtained by Milne \cite[Theorem 5.3]{mi}. Moreover,
 as noted in \cite{kr}, another equivalent formulation is 
 a special case of a 
 discrete Selberg integral
conjectured by Askey \cite{a} and proved by Evans \cite{e}. 
However, the form given by Krattenthaler is more useful for our purposes.

In  the special case
 $x=y=-q$,  \eqref{krs} evaluates the Schur $Q$-polynomial
 $Q_\lambda(1,q,\dots,q^n)$, where $\lambda$ is the odd staircase
 partition $(2m-1,2m-3,\dots,3,1)$. 
To see this, we let $x_i\equiv q^{2i-1}$ in \eqref{qse}. 
The resulting sum contains the Vandermonde determinant
$$\det_{1\leq i,j\leq m}\left(q^{(2i-1)k_j}\right)=\prod_{i=1}^mq^{2k_i}\prod_{1\leq i<j\leq m}(q^{2k_j}-q^{2k_i}) $$
and is thus  evaluated by \eqref{krs} as
\begin{align}&P_n(q,q^3,\dots,q^{2m-1})\notag\\
&\qquad\quad=
\sum_{0\leq k_m<\dots<k_1\leq n}\,\prod_{1\leq i<j\leq
  m}(q^{k_j}-q^{k_i})^2
\prod_{i=1}^m\frac{(-q;q)_{k_i}(-q;q)_{n-k_i}}{(q;q)_{k_i}(q;q)_{n-k_i}}\,(-q)^{k_i}\notag\\
&\qquad\quad=(-1)^{\binom m2}q^{\frac 14\binom{2m}3}\prod_{i=1}^m\frac{(-q,-q,q;q)_{i-1}(q^{i+m};q)_{n+1-m}}{(q;q)_{n+i-m}}\label{osq}\\
&\qquad\quad=(-1)^{\binom m2}q^{\frac 14\binom{2m}3}\prod_{i=1}^m\frac{(q^{n+1-m+i};q)_{m}}{(q;q^2)_{i-1}(q;q^2)_{i}}\notag.
\end{align}
To check the equality of the last two expressions, the elementary identity
\begin{equation}\label{epi}\prod_{i=1}^m(aq^i;q)_m=
\prod_{i=1}^m(aq^i;q)_i(aq^i;q)_{i-1},\end{equation}
with $a=1$, is useful.

Since, by anti-symmetry,
$$P_n(q,q^3,\dots,q^{2m-1})=(-1)^{\binom m2}P_n(q^{2m-1},\dots,q^3,q), $$
we obtain the following Schur $Q$-polynomial identity. Although it is
equivalent to a special case of \eqref{krs}, we have not found it
explicitly in the literature.

\begin{proposition}\label{osp}
When 
$\lambda=(2m-1,2m-3,\dots,3,1)$, then
$$Q_\lambda(1,q,\dots,q^n)=
2^mq^{\frac14\binom{2m}{3}}\prod_{i=1}^m\frac{(q^{n+1-m+i};q)_{m}}{(q;q^2)_{i-1}(q;q^2)_{i}}.
$$
\end{proposition}

We will recover Proposition \ref{osp}  twice below, see
Remark \ref{oscr} and Remark \ref{sr}. A very similar result holds for the even staircase, see Corollary \ref{escr}.

We 
write down the limit case $q\rightarrow 1$ explicitly,  starting from  \eqref{osq} with  $n$ replaced by $n-1$.

\begin{corollary}\label{osc}
When 
$\lambda=(2m-1,2m-3,\dots,3,1)$, then $Q_\lambda(1^n)$, that is, the number
of
marked shifted tableaux of shape $S(\lambda)$, with labels from the
finite alphabet $1'<1<2'<2<\dots<n'<n$, equals 
$$
2^{m^2}\prod_{i=1}^m\frac{(n+i-1)!(i-1)!}{(n+i-m-1)!(i+m-1)!}.
$$
\end{corollary}

 \begin{problem}
Give a combinatorial proof of Corollary \ref{osc}, or more generally of Proposition \ref{osp}.
\end{problem}

\section{Relation to Christoffel--Darboux kernels}
\label{ots}

The following fundamental result identifies 
 $P_n(x,y)$ with a Christoffel--Darboux kernel 
 for continuous $q$-Jacobi polynomials.

\begin{lemma}\label{tvp}
One has
\begin{equation}\label{tvpe}
P_n(x,y)=\frac{1-q^{n+1}}{1+q^{n+1}}
\frac{yP_{n+1}(x)P_{n-1}(y)-xP_{n-1}(x)P_{n+1}(y)}{1-xy}.
\end{equation}
Equivalently, 
fixing square roots $\sqrt{-x}$ and $\sqrt{-y}$
and writing
 $\xi=\sqrt{-x}+(\sqrt{-x})^{-1}$, $\eta=\sqrt{-y}+(\sqrt{-y})^{-1}$,
one has
\begin{equation}\label{pk} P_n(x,y)=\left(\sqrt {-x}\right)^{n-1}\left(\sqrt {-y}\right)^{n-1}(y-x)\,\tilde K^n(\xi,\eta),
\end{equation}
where $\tilde K_n$ is as in \eqref{ecd2}, with 
 $p_n=c_n$  defined in \eqref{mup}.
\end{lemma}

\begin{proof}
 By definition,  $(1-xy)P_n(x,y)$ equals 
\begin{multline*}(1-xy)\frac{(-q;q)_n^2}{(q;q)_n^2}\sum_{k,l=0}^n\frac{q^l-q^k}{q^l+q^k}\frac{(-q,q^{-n};q)_k}{(q,-q^{-n};q)_k}\frac{(-q,q^{-n};q)_l}{(q,-q^{-n};q)_l}\,x^ky^l\\
=\frac{(-q;q)_n^2}{(q;q)_n^2}\sum_{k,l=0}^{n+1}\frac{q^l-q^k}{q^l+q^k}\frac{(-q,q^{-n};q)_{k-1}}{(q,-q^{-n};q)_k}\frac{(-q,q^{-n};q)_{l-1}}{(q,-q^{-n};q)_l}\,x^ky^l\\
\times\Big\{
(1+q^k)(1-q^{k-n-1})(1+q^l)(1-q^{l-n-1})\\
-(1-q^k)(1+q^{k-n-1})(1-q^l)(1+q^{l-n-1})
\Big\}.
 \end{multline*}
We observe that 
the quantity within brackets factors as
$$2(q^k+q^l)(1-q^{-n-1})(1-q^{k+l-n-1}), $$
and then split the summand in a different way, writing
\begin{multline*}2(q^l-q^k)(1+q^{-n-1})(1-q^{k+l-n-1})\\
=(1-q^k)(1-q^{k-n-1})(1+q^l)(1+q^{l-n-1})-(1+q^k)(1+q^{k-n-1})(1-q^l)(1-q^{l-n-1}).
\end{multline*}
This gives
\begin{multline*}(1-xy)P_n(x,y)\\
=\frac{(-q;q)_n^2}{(q;q)_n^2}\left(
\sum_{k=1}^n\frac{(-q;q)_{k-1}(q^{-n};q)_k}{(q;q)_{k-1}(-q^{-n};q)_k}\,x^k\sum_{l=0}^{n+1}\frac{(-q,q^{-n-1};q)_{l}}{(q,-q^{-n-1};q)_l}\,y^l\right.\\
\left.-\sum_{k=0}^{n+1}\frac{(-q,q^{-n-1};q)_{k}}{(q,-q^{-n-1};q)_{k}}\,x^k\sum_{l=1}^n\frac{(-q;q)_{l-1}(q^{-n};q)_{l}}{(q;q)_{l-1}(-q^{-n};q)_{l}}\,y^l\right).
\end{multline*}
Replacing $k$ by $k+1$ in the first term and $l$ by $l+1$ in the
second term yields \eqref{tvpe}. It is then straight-forward to derive \eqref{pk}, using \eqref{pmu}, \eqref{cn} and the elementary identity
$$\xi^2-\eta^2=\frac 1{xy}(x-y)(1-xy).$$
\end{proof}

\begin{remark}\label{kr}
By the $q$-binomial theorem \cite[(II.3)]{gr},
$$\lim_{n\rightarrow\infty}P_n(x)=\frac{(-q;q)_\infty}{(q;q)_\infty}\sum_{k=0}^\infty\frac{(-q;q)_k}{(q;q)_k}(-x)^k=
\frac{(-q,qx;q)_\infty}{(q,-x;q)_\infty},
 $$
for $|q|,|x|<1$. Lemma \ref{tvp} then gives
$$\lim_{n\rightarrow \infty}P_n(x,y)=\frac{(-q;q)_\infty^2}{(q;q)_\infty^2}\frac{(qx,qy;q)_\infty}{(-x,-y;q)_\infty}\frac{y-x}{1-xy}. $$
By  \eqref{pfi} and \eqref{spb},
it follows that
\begin{equation}\label{k1}\lim_{n\rightarrow\infty}P_n(x_1,\dots,x_{m})=\frac{(-q;q)_\infty^m}{(q;q)_\infty^m}\prod_{i=1}^{m}\frac{(qx_i;q)_\infty}{(-x_i;q)_\infty}\prod_{1\leq
  i<j\leq {m}}\frac{x_j-x_i}{1-x_ix_j} \end{equation}
for $m$ even.  By \eqref{peo}, this holds also for odd $m$.
Comparing \eqref{k1} and \eqref{k2} yields
 Kawanaka's identity \eqref{kbf}, which thus follows from  Lemma \ref{tvp}. In fact, our proof of  Lemma \ref{tvp} generalizes
   the proof of Kawanaka's identity given by Ishikawa and Wakayama
  \cite{iw2}. 
\end{remark}

We rewrite \eqref{ktk} in terms of the polynomials $P_n$, assuming for convenience that we have chosen $\sqrt x$ and $\sqrt y$ so that $\sqrt x/\sqrt y=\sqrt{-x}/\sqrt{-y}$.  The resulting identity will be used in the proof of Theorem \ref{th}.

\begin{corollary}\label{rce}
One has
$$P_n(x,y)=\frac{\sqrt x+\sqrt
  y}{1+\sqrt{xy}}
\sum_{j=0}^{n-1}(\sqrt{xy})^{n-1-j}
\big(\sqrt y P_{j+1}(x)P_j(y)-\sqrt x P_j(x)P_{j+1}(y)\big).$$
\end{corollary}

Next, combining \eqref{pk} and \eqref{pfi} gives 
$$P_n(x_1,\dots,x_{2m})=
\prod_{i=1}^{2m}(\sqrt{-x_i})^{n-1}\pfaff_{1\leq i,j\leq 2m}((x_j-x_i)K^n(\xi_i,\xi_j)),
 $$
where $\xi_i=\sqrt{-x_i}+(\sqrt{-x_i})^{-1}$. We see that 
 the  pfaffian is of the form
 \eqref{kpb}, with $\zeta_i=x_i$, $z_i=\xi_i$. This shows the following result when $m$ is even. The case off odd  $m$  follows as a limit case, using \eqref{peo} and \eqref{sfl}, and the expression \eqref{cn} for the norms.

\begin{theorem}\label{pkt}
For any fixed choice of square roots $\sqrt{-x_i}$, 
\begin{multline*}P_n(x_1,\dots,x_m)\\
=
\prod_{i=1}^{m}\frac{(-q;q)_{n+1-i}}{(q;q)_{n+1-i}}(\sqrt{-x_i})^{n+1-m}
\prod_{1\leq i<j\leq m}(x_j-x_i)\,
 \tilde{\mathbb  P}_{(n+1-m)^m}(\xi_1,\dots,\xi_{m}),
 \end{multline*}
where $\xi_i=\sqrt{-x_i}+(\sqrt{-x_i})^{-1}$, and where  $\tilde{\mathbb  P}$ is defined by choosing  $p_k=c_k$ in \eqref{mrs}. When $m$ is even, 
we equivalently have
$$P_n(x_1,\dots,x_{2m})
=\prod_{i=1}^{2m}(\sqrt{-x_i})^{n+1-2m}
\prod_{1\leq i<j\leq 2m}(x_j-x_i)\,
 \tilde K_m^n(\xi_1,\dots,\xi_{2m}),
$$
where $\tilde K$ is defined by choosing $p_k=c_k$ in \emph{Corollary \ref{cdc}}.
\end{theorem}

We will give an alternative proof of Theorem \ref{pkt} in the Appendix.

Note that $\xi_i$ is invariant under the transformations $\sqrt{-x_i}\mapsto(\sqrt{-x_i})^{-1} $ or, equivalently, $\sqrt{x_i}\mapsto-(\sqrt{x_i})^{-1}$. 
This yields the following  hyperoctahedral symmetry, which is  quite non-obvious from the definition \eqref{pd}.

\begin{corollary}\label{hoc}
The function
$$\prod_{i=1}^m(\sqrt{x_i})^{m-n-1}\frac{P_n(x_1,\dots,x_m)}{\prod_{1\leq i<j\leq m}(x_j-x_i)} $$
is invariant under the action of the hyperoctahedral group $B_m$ generated by permutation of the variables together with inversions $\sqrt{x_i}\mapsto -(\sqrt{x_i})^{-1}$.
 \end{corollary}

It should be noted that, although phrased in terms of square roots, this is a  polynomial statement. For instance, the symmetry under $\sqrt{x_1}\mapsto -(\sqrt{x_1})^{-1}$ can be rephrased as
$$\prod_{i=2}^m(1-x_1x_i)\,P_n(x_1,x_2,\dots,x_m)=(-x_1)^n\prod_{i=2}^m(x_i-x_1)\,P_n(x_1^{-1},x_2,\dots,x_m).$$

We note in passing the  formulas for $P_n$ obtained from \eqref{mk}.
These will not be used in the remainder of this paper.

\begin{corollary} One has
\begin{multline*}P_n(x_1,\dots,x_m,y_1,\dots,y_m)=\frac{\prod_{i,j=1}^m(y_j-x_i)}{\prod_{1\leq i<j\leq m}(1-x_ix_j)(1-y_iy_j)}\,\det_{1\leq i,j\leq m}\left(\frac{P_n(x_i,y_j)}{y_j-x_i}\right)\\
=\frac{\prod_{i,j=1}^m(y_j-x_i)}{\prod_{1\leq i<j\leq m}(1-x_ix_j)(1-y_iy_j)}
\,\sum_{\substack{0\leq k_m<\dots<k_1\leq n-1\\k_1,\dots,k_m\equiv n-1\, \operatorname{mod}\, 2}}\Bigg(\,\prod_{i=1}^m\frac{1+q^{k_i+1}}{1-q^{k_i+1}}\\
\times\det_{1\leq i,j\leq m}\left(x_j^{(n-1-k_i)/2}P_{k_i}(x_j)\right)\det_{1\leq i,j\leq m}\left(y_j^{(n-1-k_i)/2}P_{k_i}(y_j)\right)\Bigg).\end{multline*}
\end{corollary}

Turning to \eqref{kpa}, we have
$$P_n(x_1,\dots,x_{2m})=\frac{\prod_{1\leq i<j\leq 2m}(x_j-x_i)}{\prod_{i=1}^{2m} x_i^{m-1}\prod_{1\leq i<j\leq 2m}(w_j-w_i)}\pfaff_{1\leq i,j\leq 2m}\left(\frac{w_j-w_i}{x_j-x_i}\,P_n(x_i,x_j)\right), $$
where
$w_i=\sqrt{x_i}-(\sqrt{x_i})^{-1}$, so that
$$w_j-w_i=\frac{1}{\sqrt{x_ix_j}}(\sqrt{x_j}-\sqrt{x_i})(1+\sqrt{x_ix_j}). $$
This gives the following pfaffian formula, which forms the basis for Theorem \ref{th}.

\begin{corollary}\label{api}
For $m$ even,
$$P_n(x_1,\dots,x_{m})=\prod_{1\leq i<j\leq m}\frac{\sqrt {x_i}+\sqrt{x_j}}{1+\sqrt{x_ix_j}}
\pfaff_{1\leq i,j\leq m}\left(\frac{1+\sqrt{x_ix_j}}{\sqrt {x_i}+\sqrt{x_j}}\,P_n(x_i,x_j)\right). $$
\end{corollary}

Rewriting the determinant formula \eqref{mrs} in terms of the one-variable polynomials $P_n(x)$ and  reversing the order of the columns, we obtain the  following identity. 

\begin{corollary}\label{dt}
One has
\begin{multline}\label{dte}
\prod_{1\leq i<j\leq m}(1-x_ix_j)\,P_n(x_1,\dots,x_m)\\
=\prod_{i=1}^m\frac{(q^{n+1-m+i};q)_{i-1}}{(-q^{n+1-m+i};q)_{i-1}}
\det_{1\leq i,j\leq m}\left(x_i^{j-1}P_{n+m+1-2j}(x_i)\right).
\end{multline}
\end{corollary}

Corollary \ref{dt} seems interesting enough to write  down also in standard notation for Schur $Q$-polynomials.

\begin{corollary}\label{qdc}
For $\lambda$ a partition of length $m$, 
\begin{multline*}
Q_\lambda(1,q,\dots,q^n)
=\prod_{i=1}^m\frac{(q^{n+1-m+i};q)_{i-1}}{(-q^{n+1-m+i};q)_{i-1}}\prod_{1\leq i<j\leq m}\frac{1}{1-q^{\lambda_i+\lambda_j}}\\
\times\det_{1\leq i,j\leq m}\left(q^{(j-1)\lambda_i}Q_{(\lambda_i)}(1,q,\dots,q^{n+m+1-2j})\right).
\end{multline*}
\end{corollary}

 Replacing $n$ by $n-1$ and letting $q\rightarrow 1$ gives the following 
simple determinant formula for $Q_\lambda(1^n)$. 

\begin{corollary}\label{fdc}
For $\lambda$ a partition of length $m$,
$$
Q_\lambda(1^n)
=2^{\frac 12m(2n+1-m)}\prod_{i=1}^m\frac{1}{(n-m+i-1)!}\prod_{1\leq i<j\leq m}\frac{1}{\lambda_i+\lambda_j}\,
\det_{1\leq i,j\leq m}\left(f_{n+m-2j}(\lambda_i))\right),$$
where $f_k$ are the hypergeometric polynomials defined in \eqref{fk}.
\end{corollary}

For future refrence \cite{rss}, 
we rewrite
Corollary \ref{fdc} in terms of  generalized Schur polynomials and  multivariable Christoffel--Darboux kernels corresponding to the continuous dual Hahn polynomials $p_k^{(0)}$ and $p_k^{(1)}$ introduced in \eqref{cdh}.
Let $\mathbb P_{n^m}^{(\ep)}$ be the polynomial obtained by choosing $p_k=p_k^{(\ep)}$ in \eqref{rsd}, and let $K_m^{n,\ep}$ be the kernel similarly obtained from Lemma \ref{cdl}. Using the expression \eqref{cdhn} for the norms, we arrive at the following result, which can  also be obtained 
directly from Theorem \ref{pkt} as a limit case.

\begin{corollary}
Let $n-m=2k+\ep$, with $k$ an integer and $\ep\in\{0,1\}$. Then,
\begin{multline*}Q_{(\lambda_1,\dots,\lambda_m)}(1^n)\\
=\frac{2^{\frac 12 m(2n+1-m)}(-1)^{km}}
{\prod_{i=1}^m{(n-i)!}}
\prod_{i=1}^m\lambda_i^\ep\prod_{1\leq i<j\leq m}(\lambda_i-\lambda_j)
\,\mathbb P^{(\ep)}_{k^m}(-\lambda_1^2,\dots,-\lambda_m^2).
 \end{multline*}
In particular,
$$Q_{(\lambda_1,\dots,\lambda_{2m})}(1^{2n})=4^m\prod_{1\leq i<j\leq 2m}(\lambda_i-\lambda_j)\,K_m^{n,0}(-\lambda_1^2,\dots,-\lambda_{2m}^2),
 $$
$$Q_{(\lambda_1,\dots,\lambda_{2m})}(1^{2n+1})=4^m\prod_{i=1}^{2m}\lambda_i\prod_{1\leq i<j\leq 2m}(\lambda_i-\lambda_j)\,K_m^{n,1}(-\lambda_1^2,\dots,-\lambda_{2m}^2).
 $$
\end{corollary}

\section{A pfaffian hypergeometric identity}\label{mss}

Using the pfaffian formula of  Corollary \ref{api}, we derive the following hypergeometric identity. Note that it displays the hyperoctahedral symmetry of Corollary \ref{hoc}.

\begin{theorem}\label{th} For any fixed choice of $\sqrt{x_i}$, 
\begin{multline}\label{the}
P_n(x_1,\dots,x_m)=q^{\frac 14 m(m-1-2n)}\left(\frac{1-q^{n+1}}{1-q}\right)^m
\prod_{i=1}^m (\sqrt{x_i})^n\prod_{1\leq i<j\leq m} 
\frac{\sqrt{x_i}+\sqrt
    {x_j}}{1+\sqrt{x_ix_j}}
 \\
\times\sum_{k_1,\dots,k_m=0}^n\,\prod_{1\leq i<j\leq m} 
\frac{q^{k_j}-q^{k_i}}{1-q^{k_i+k_j+1}}\prod_{i=1}^m\frac{(q^{-n},q^{n+2},\sqrt{q/x_i},-\sqrt{qx_i};q)_{k_i}}{(q,q,q^{3/2},-q^{3/2};q)_{k_i}}\,q^{k_i}.
\end{multline}
\end{theorem}

\begin{proof}

We first consider the case $m=2$. Using \eqref{43} in Corollary \ref{rce}
and changing the order of summation gives
\begin{multline}\label{pts}
P_n(x,y)=\frac{q^{-\frac 12}}{(1-q)^2}\,(\sqrt{xy})^{n}
\frac{\sqrt x+\sqrt y}{1+\sqrt {xy}}\\
\times\sum_{k,l=0}^n\frac{(\sqrt {q/x},-\sqrt{qx};q)_k}{(q,q,q^{3/2},-q^{3/2};q)_k}\frac{(\sqrt{q/y},-\sqrt{qy};q)_l}{(q,q,q^{3/2},-q^{3/2};q)_l}\,q^{k+l}\\
\times\sum_{j=\max(k-1,l-1)}^{n-1}(1-q^{j+1})(1-q^{j+2})\,q^{-j}\\
\times\big\{(q^{-j-1},q^{j+3};q)_k(q^{-j},q^{j+2};q)_l-
(q^{-j},q^{j+2};q)_k(q^{-j-1},q^{j+3};q)_l
\big\}.
\end{multline}

Let $S$ denote the inner sum in \eqref{pts}. It can be rewritten as
\begin{multline*}S=\sum_{j=\max(k-1,l-1)}^{n-1}(-1)^{k+l}q^{\binom k2+\binom
  l2-(k+l)(j+1)-j}
(q^{j-k+2};q)_{2k}(q^{j-l+2};q)_{2l}\\
\times\big\{
q^l(1-q^{j+k+2})(1-q^{j-l+1})-q^k(1-q^{j-k+1})(1-q^{j+l+2})
\big\},
\end{multline*}
where the expression within brackets factors as
$$(q^l-q^k)(1-q^{2j+3}).$$
Replacing $j$ by $n-1-j$, elementary
manipulations give
\begin{multline*}
S=q^{1-n}(q^l-q^k)(1-q^{2n+1})(q^{-n},q^{n+1};q)_k(q^{-n},q^{n+1};q)_l\\
\times\sum_{j=0}^{\min(n-k,n-l)}\frac{1-q^{-2n-1+2j}}{1-q^{-2n-1}}\frac{(q^{k-n},q^{l-n};q)_j}{(q^{-k-n},q^{-l-n};q)_j}\,q^{-(k+l+1)j}.
\end{multline*}
Here, the sum is a very-well-poised ${}_6\phi_5$, which by \cite[Eq.\
(II.21)]{gr} equals
$$\frac{(1-q^{n+k+1})(1-q^{n+l+1})}{(1-q^{2n+1})(1-q^{k+l+1})},$$
so that
$$S=q^{1-n}(1-q^{n+1})^2\frac{q^l-q^k}{1-q^{k+l+1}}\,(q^{-n},q^{n+2};q)_k(q^{-n},q^{n+2};q)_l.
$$
This gives
\begin{multline}\label{nds}
P_n(x,y)=q^{\frac 12-n}\left(\frac{1-q^{n+1}}{1-q}\right)^2(\sqrt{xy})^{n}
\frac{\sqrt x+\sqrt y}{1+\sqrt
  {xy}}\\
\times\sum_{k,l=0}^n\left(\frac{q^l-q^k}{1-q^{k+l+1}}
\frac{(q^{-n},q^{n+2},\sqrt{q/x},-\sqrt{qx};q)_k}{(q,q,q^{3/2},-q^{3/2};q)_k}
\,q^k\right.\\
\times\left.\frac{(q^{-n},q^{n+2},\sqrt{q/y},-\sqrt{qy};q)_l}{(q,q,q^{3/2},-q^{3/2};q)_l}\,q^{l}\right),
\end{multline}
which is the case $m=2$ of \eqref{the}.

Suppose now that $m$ is even. We  plug \eqref{nds} into Corollary \ref{api}, thus expressing $P_n(x_1,\dots,x_m)$ as a pfaffian of two-dimensional sums. 
Applying Corollary \ref{lms} leads to an $m$-dimensional sum, each term containing a pfaffian of the form \eqref{spb}.
This proves \eqref{the} in the case when $m$ is even.

Finally, if $m$ is odd, we use \eqref{peo} to write
$$P_n(x_1,\dots,x_m)=\lim_{t\rightarrow 0}\frac{(q;q)_{n+1}}{(-q;q)_{n+1}}\frac {(-1)^m}{x_1\dotsm x_m}\,P_{n+1}(x_1,\dots,x_m,t). $$
Applying \eqref{the} to the right-hand side, only the term with $k_{m+1}=n+1$ survives in the limit. After simplification, one arrives at 
 the desired expression.
\end{proof}

In contrast to \eqref{pd}, \eqref{the} is regular in the limit
\begin{equation}\label{lim}Q_\lambda(1^{n+1})=2^m\lim_{q\rightarrow 1}
P_n(q^{\lambda_1},\dots,q^{\lambda_m}).\end{equation}

\begin{corollary}\label{qc}
Let $\lambda$ be a partition  of length $m$. Then,
$$Q_\lambda(1^{n+1})=
(2n+2)^m\sum_{k_1,\dots,k_m=0}^n\,\prod_{1\leq i<j\leq m}\frac{k_i-k_j}{k_i+k_j+1}\prod_{i=1}^m\frac{(-n,n+2,(1-\lambda_i)/2)_{k_i}}{(1,1,3/2)_{k_i}}.
$$
\end{corollary}

\begin{remark}\label{oscr} If $x_i\equiv q^{2i-1}$, corresponding to the  odd staircase partition, the sum in \eqref{the} reduces to the term with $k_i\equiv i-1$. After simplification, one recovers
 Proposition \ref{osp}. 
\end{remark}

\section{Determinantal hypergeometric identities}
\label{ds}

By inserting hypergeometric expressions  for $P_n(x)$ into the determinant in 
Corollary~\ref{dt}, one may obtain further formulas for $P_n(x_1,\dots,x_m)$.
We are particularly interested in  cases when the result simplifies using a determinant evaluation such as
\cite[Lemma A.1]{sc}
\begin{multline}\label{sd}
\det_{1\leq i,j\leq m}\left(\frac{(AX_i,AC/X_i;q)_{j-1}}{(BX_i,BC/X_i;q)_{j-1}}\right)=q^{\binom m3}(AC)^{\binom m2}\\
\times\prod_{1\leq i<j\leq m}(X_j-X_i)(1-X_iX_j/C)
\prod_{i=1}^m\frac{(B/A,ABCq^{2m-2i};q)_{i-1}}{X_i^{m-1}(BX_i,BC/X_i;q)_{m-1}}
\end{multline}
or the degenerate case
\begin{multline}\label{sdd}\det_{1\leq i,j\leq m}\left((AX_i;q)_{j-1}(BX_i;q)_{m-j}\right)\\
=q^{\binom m3}A^{\binom m2}\prod_{1\leq i<j\leq m}(X_i-X_j)
\prod_{i=1}^m(q^{i-m}B/A;q)_{i-1}.
\end{multline}
Multiple hypergeometric sums obtained in this way have been called ``Schlosser-type'', see \cite{gk,sc}.

First of all, using our original definition of $P_n(x)$, we  find 
the following identity, which has a very similar structure to sums found in \cite{rs}.

\begin{theorem}\label{ot}
One has
\begin{multline}\label{ote}
\prod_{1\leq i<j\leq m}(1-x_ix_j)P_n(x_1,\dots,x_m)\\
\begin{split}&=
(-1)^{\binom m2}q^{(n+1)\binom m2+2\binom m3}\frac{(-q;q)_n^m}{(q;q)_{m+n-1}^m}\prod_{i=1}^m\frac{(q^{n+1-m+i};q)_{i-1}}{(-q^{2-m};q)_{i-1}}\\
&\quad\times\sum_{k_1,\dots,k_m=0}^{n+m-1}\left(
\prod_{1\leq i<j\leq m}(q^{k_j}-q^{k_i})(1-q^{k_i+k_j+1-m-n})\right.\end{split}\\
\times\left.\prod_{i=1}^m\frac{(q^{1-m-n},-q^{2-m};q)_{k_i}}{(q,-q^{-n};q)_{k_i}}
\, x_i^{k_i}\right).
\end{multline}
\end{theorem}

\begin{proof}
By \eqref{ovp}, the  right-hand side of \eqref{dte} equals
$$\prod_{i=1}^m\frac{(-q;q)_{n+1-i}}{(q;q)_{n+1-i}}
\sum_{\sigma\in S_m}\sgn(\sigma) \sum_{k_1,\dots,k_m\geq 0}
\prod_{i=1}^m\frac{(-q,q^{2\sigma(i)-n-m-1};q)_{k_i}}{(q,-q^{2\sigma(i)-n-m-1};q)_{k_i}}\,x_i^{k_i+\sigma(i)-1}.
 $$
Replacing $k_i$ by $k_i+1-\sigma(i)$ and interchanging summations gives after simplification
\begin{multline*}(-1)^{\binom m2}\frac{(-q;q)_{m+n-1}^m}{(q;q)_{m+n-1}^m}\prod_{i=1}^m
\frac{(q^{n+1-m+i};q)_{i-1}}{(-q^{n+1-m+i};q)_{i-1}}\\
\times\sum_{k_1,\dots,k_m=0}^{m+n-1}\det_{1\leq i,j\leq m}\left(\frac{(q^{k_i+1-m-n},q^{-k_i};q)_{j-1}}{(-q^{k_i+1-m-n},-q^{-k_i};q)_{j-1}}\right)\prod_{i=1}^m\frac{(-q,q^{1-m-n};q)_{k_i}}{(q,-q^{1-m-n};q)_{k_i}}\,x_i^{k_i}.
 \end{multline*}
By \eqref{sd}, the determinant equals
$$q^{\binom m3}\prod_{1\leq i<j\leq m}(q^{k_j}-q^{k_i})(1-q^{k_i+k_j+1-m-n})
\prod_{i=1}^m\frac{(-1,-q^{1+m-n-2i};q)_{i-1}}{(-q^{k_i+1-m-n},-q^{-k_i};q)_{m-1}}\,q^{(1-m)k_i}. $$
After further simplification, one arrives at the right-hand side of \eqref{ote}.\end{proof}

\begin{remark}
Writing Theorem \ref{ot} as
$$\prod_{1\leq i<j\leq m}(1-x_ix_j)\,P_n(x_1,\dots,x_m)
=\sum_{k_1,\dots,k_m=0}^{n+m-1}\chi(k_1,\dots,k_m)\prod_{i=1}^m x_i^{k_i},$$
one easily checks that 
$(-1)^{|k|}\chi(k)$ is antisymmetric and invariant under the reflections $k_i\mapsto n+m-1-k_i$. This fact is  equivalent to Corollary \ref{hoc}. 
In analogy with the  Schur polynomial expansion \eqref{pse}, one may exploit this symmetry to rewrite Theorem \ref{ot}  in terms of characters of the classical groups $\mathrm{SO}(2m)$ and $\mathrm{SO}(2m+1)$, according to whether $n+m$  is odd or even, respectively.
\end{remark}

Inserting different expressions for $P_n(x)$ in Corollary \ref{dt}, one may  obtain further hypergeometric formulas. 
We will not  exploit all possibilities, but merely give a few examples,
corresponding to the eight expressions in \eqref{chia}. The resulting identities  share the nice properties of  Theorem
 \ref{th}, that is, they
are 
regular in the limit $q\rightarrow 1$, have completely factored terms, display hyperoctahedral symmetry, and have nice specializations at staircase-type partitions.

For instance, assuming that $n+m$ is odd, we  insert \eqref{chiaa} into \eqref{dte}. After simplification, we obtain
\begin{multline*}\prod_{1\leq i<j\leq m}{(1-x_ix_j)}P_n(x_1,\dots,x_m) \\
=\frac{q^{-\frac 12 nm}}{(1-q)^m}
\prod_{i=1}^m\frac{(q^{n+1-m+i};q)_i}{(-q^{n+1-m+i};q)_{i-1}(q^{(1-n-m)/2},-q^{(3+n-m)/2};q)_{i-1}}\, x_i^{(n+m-1)/2}\\
\times
\sum_{k_1,\dots,k_m=0}^{(n+m-1)/2}\,
\prod_{i=1}^m\frac{(q^{(1-n-m)/2},-q^{(3+n-m)/2},qx_i,q/x_i;q)_{k_i}}
{(q,q,q^{3/2},-q^{3/2};q)_{k_i}}\,q^{k_i}\\
\times\det_{1\leq i,j\leq m}\left((q^{k_i+(1-n-m)/2};q)_{j-1}
(-q^{k_i+(3+n-m)/2};q)_{m-j}
\right),
\end{multline*}
where, by \eqref{sdd}, the determinant equals
$$q^{-\frac 12\binom m3-\frac{n+1}{2}\binom m2}\prod_{1\leq i<j\leq m}(q^{k_i}-q^{k_j})
\prod_{i=1}^m(-q^{n+1-m+i};q)_{i-1}.
 $$
Repeating this for all eight expressions  \eqref{chia} yields the following result. In the case of \eqref{ec} and \eqref{od} we have also used \eqref{epi}.

\begin{theorem}\label{dft}
If $n+m$ is odd, then $\prod_{1\leq i<j\leq m}{(1-x_ix_j)}P_n(x_1,\dots,x_m)$ is equal to each of the four expressions
\begin{subequations}
\begin{multline}
\frac{(-1)^{\binom m2}q^{-\frac 12\binom{m+1}3-\frac n2\binom{m+1}2}}{(1-q)^m}
\prod_{i=1}^m\frac{(q^{n+1-m+i};q)_i}{(q^{(1-n-m)/2},-q^{(3+n-m)/2};q)_{i-1}}\, x_i^{(n+m-1)/2}\\
\times
\sum_{k_1,\dots,k_m=0}^{(n+m-1)/2}\left(\prod_{1\leq i<j\leq m}(q^{k_j}-q^{k_i})\right.\\
\times\left.\prod_{i=1}^m\frac{(q^{(1-n-m)/2},-q^{(3+n-m)/2},qx_i,q/x_i;q)_{k_i}}
{(q,q,q^{3/2},-q^{3/2};q)_{k_i}}\,q^{k_i}\right),\label{ea}
\end{multline}
\begin{multline}
\frac{q^{-\frac 14 nm}(-q;q)_{(n+1-m)/2}^m}{(q;q)_{(n+m-1)/2}^m}
\prod_{i=1}^m\frac{(q^{n+1-m+i};q)_{i-1}(q^{3/2};q)_{(n+m+1-2i)/2}}{(-q^{1/2};q)_{(n+m+1-2i)/2}}\, x_i^{(n+m-1)/2}\\
\times\sum_{k_1,\dots,k_m=0}^{(n+m-1)/2}\left(\prod_{1\leq i<j\leq m}(q^{k_j}-q^{k_i})\right.\\
\times\left.\prod_{i=1}^m\frac{(q^{(1-n-m)/2},-q^{(3+n-m)/2},q^{1/2}x_i,q^{1/2}/x_i;q)_{k_i}}
{(q,q^{1/2},q^{3/2},-q;q)_{k_i}}\,q^{k_i}\right),\label{eb}
\end{multline}
\begin{multline}
\frac{q^{-\frac 12 nm}}{(1-q)^m(q^{3+n-m};q^2)_{m-1}^m}
\prod_{i=1}^m(q^{n+1-m+i};q)_{m}\, x_i^{(n+m-1)/2}\\
\times\sum_{k_1,\dots,k_m=0}^{(n+m-1)/2}\left(\prod_{1\leq i<j\leq m}(q^{2k_j}-q^{2k_i})\right.
\\
\times\left.\prod_{i=1}^m\frac{(q^{1-n-m},q^{3+n-m},qx_i,q/x_i;q^2)_{k_i}}
{(q,q^2,q^2,q^3;q^2)_{k_i}}\,q^{2k_i}\right),\label{ec}
\end{multline}
\begin{multline}
\frac{q^{- nm}}{(1-q)^{2m}(q^{3+n-m};q^2)_{m-1}^m}
\prod_{i=1}^m(q^{n+1-m+i};q)_{i}^2\, x_i^{(n+m-1)/2}
\\
\times\sum_{k_1,\dots,k_m=0}^{(n+m-1)/2}\left(\prod_{1\leq i<j\leq m}(q^{2k_j}-q^{2k_i})\right.
\\
\times\left.\prod_{i=1}^m\frac{(q^{1-n-m},q^{3+n-m},q^2x_i,q^2/x_i;q^2)_{k_i}}
{(q^2,q^2,q^3,q^3;q^2)_{k_i}}\,q^{2k_i}\right),\label{ed}
\end{multline}
\end{subequations}
whereas if $n+m$ is even, $\prod_{1\leq i<j\leq m}{(1-x_ix_j)}P_n(x_1,\dots,x_m)$ equals
\begin{subequations}
\begin{multline}
\frac{(-1)^{\binom m2}q^{-\frac 12\binom{m+1}3-\frac {n-1}2\binom{m+1}2}}{(1-q)^{2m}}
\prod_{i=1}^m\frac{(q^{n+1-m+i};q)_i}{(q^{(2-n-m)/2},-q^{(4+n-m)/2};q)_{i-1}}\\
\times\prod_{i=1}^m(1-x_i) x_i^{(n+m-2)/2}
\sum_{k_1,\dots,k_m=0}^{(n+m-2)/2}\left(\prod_{1\leq i<j\leq m}(q^{k_j}-q^{k_i})\right.\\
\times\left.\prod_{i=1}^m\frac{(q^{(2-n-m)/2},-q^{(4+n-m)/2},qx_i,q/x_i;q)_{k_i}}
{(q,q^2,q^{3/2},-q^{3/2};q)_{k_i}}\,q^{k_i}\right),
\label{oa}\end{multline}
\begin{multline}\label{ob}
\frac{q^{-\frac 14 (n-1)m}(-q;q)_{(n+2-m)/2}^m}{(1-q)^m(q;q)_{(n+m-2)/2}^m}
\prod_{i=1}^m\frac{(q^{n+1-m+i};q)_{i-1}(q^{3/2};q)_{(n+m-2i)/2}}{(-q^{3/2};q)_{(n+m-2i)/2}}\\
\times\prod_{i=1}^m(1-x_i) x_i^{(n+m-2)/2}\sum_{k_1,\dots,k_m=0}^{(n+m-2)/2}\left(\prod_{1\leq i<j\leq m}(q^{k_j}-q^{k_i})\right.\\
\times\left.\prod_{i=1}^m\frac{(q^{(2-n-m)/2},-q^{(4+n-m)/2},q^{1/2}x_i,q^{1/2}/x_i;q)_{k_i}}
{(q,q^{3/2},q^{3/2},-q;q)_{k_i}}\,q^{k_i}\right),
\end{multline}
\begin{multline}\label{oc}
\frac{q^{-\frac 12 (n-1)m}}{(1-q)^{2m}(q^{4+n-m};q^2)_{m-2}^m}
\prod_{i=1}^m(q^{n+1-m+i};q)_{i-1}^2\\
\times\prod_{i=1}^m(1-x_i) x_i^{(n+m-2)/2}\sum_{k_1,\dots,k_m=0}^{(n+m-2)/2}\left(\prod_{1\leq i<j\leq m}(q^{2k_j}-q^{2k_i})\right.
\\
\times\left.\prod_{i=1}^m\frac{(q^{2-n-m},q^{4+n-m},qx_i,q/x_i;q^2)_{k_i}}
{(q^2,q^2,q^3,q^3;q^2)_{k_i}}\,q^{2k_i}\right),
\end{multline}
\begin{multline}\label{od}
\frac{q^{- (n-1)m}}{(1-q)^{2m}(1-q^2)^m(q^{4+n-m};q^2)_{m-2}^m}
\prod_{i=1}^m(q^{n+1-m+i};q)_{m}
\\
\times\prod_{i=1}^m(1-x_i) x_i^{(n+m-2)/2}\sum_{k_1,\dots,k_m=0}^{(n+m-2)/2}\left(\prod_{1\leq i<j\leq m}(q^{2k_j}-q^{2k_i})\right.
\\
\times\left.\prod_{i=1}^m\frac{(q^{2-n-m},q^{4+n-m},q^2x_i,q^2/x_i;q^2)_{k_i}}
{(q^2,q^3,q^3,q^4;q^2)_{k_i}}\,q^{2k_i}\right).
\end{multline}\end{subequations}
\end{theorem}

Note that  all eight expressions are regular in the limit \eqref{lim}. We do not write down the resulting
 identities for $Q_\lambda(1^n)$ explicitly.

\begin{remark}\label{sr}
By specializing the variables $x_i$, each of the eight sums in Theorem~\ref{dft} may be reduced to the term with $k_i\equiv i-1$. In the case of \eqref{ea} and \eqref{oa}, this happens if $x_i\equiv q^i$, leading to the identity
\begin{equation*}\begin{split}Q_{(m,m-1,\dots,1)}(1,q,\dots,q^n)&=2^m(-1)^{\binom m2}P_n(q,q^2,\dots,q^m)\\
&=
2^mq^{\binom{m+1}3}\prod_{i=1}^m\frac{(q^{n+1-m+i};q)_i}{(q;q^2)_i},\end{split}\end{equation*}
which follows more easily from 
\cite[Ex.\ III.8.3]{m}
$$Q_{(m,m-1,\dots,1)}=2^m s_{(m,m-1,\dots,1)}.$$
In the case of \eqref{eb} and \eqref{ob}, one should let
$x_i\equiv q^{i-\frac 12}$. We write down the result in Corollary \ref{hsc}. 
Similarly, 
when $x_i\equiv q^{2i-1}$,   \eqref{ec} and \eqref{oc} reduce to
 Proposition~\ref{osp} and when $x_i\equiv q^{2i}$, \eqref{ed} and \eqref{od} yield Corollary~\ref{escr}.
\end{remark}

\begin{corollary}\label{hsc}
If  $x_i\equiv q^{i-\frac 12}$, then $P_n(x_1,\dots,x_m)$ equals
$$(-1)^{\binom m2}q^{\frac 18\binom{2m}3}\prod_{i=1}^m\frac{(q,q^{n+1-m+i};q)_{i-1}(q^{3/2},-q;q)_{(n+m+1-2i)/2}}{(-q,q^{1/2},q^{3/2};q)_{i-1}(q,-q^{1/2};q)_{(n+m+1-2i)/2}} $$
if  $n+m$ is odd and
$$(-1)^{\binom m2}q^{\frac 18\binom{2m}3}\frac{(q^{1/2};q)_m}{(1-q)^m}
\prod_{i=1}^m\frac{(q,q^{n+1-m+i};q)_{i-1}(q^{3/2};q)_{(n+m-2i)/2}(-q;q)_{(n+m+2-2i)/2}}{(-q,q^{3/2},q^{3/2};q)_{i-1}(-q^{3/2},q;q)_{(n+m-2i)/2}} $$
if $n+m$ is even.
\end{corollary}

\begin{corollary}\label{escr} When $\lambda=(2m,2m-2,\dotsm, 2)$, then
$$Q_\lambda(1,q,\dots,q^n)=2^mq^{2\binom{m+1}3 }\prod_{i=1}^m\frac{(q^{n+1-m+i};q)_i^2}{(q;q^2)_i^2}. $$
\end{corollary}

\section*{Appendix. Alternative proof of Theorem \ref{pkt}}

\setcounter{section}{1}
\renewcommand{\thesection}{\Alph{section}}
\setcounter{equation}{0}
\setcounter{theo+}{0}

The proof of Theorem \ref{pkt} given above used the explicit expression for $q$-ultra\-spherical polynomials as basic hypergeometric series. Here, we  give a different proof using only their orthogonality measure. This might be useful for the purpose of generalization.

Let
$$w(\theta)=\frac{(q;q)_\infty^2(e^{2i\theta},e^{-2i\theta};q)_\infty}{2\pi(-q;q)_\infty^2(-qe^{2i\theta},-qe^{-2i\theta};q)_\infty}. $$
By \cite[Eq.\ (7.4.15)]{gr},
\begin{equation}\label{oi}\int f(x)\,d\mu(x)=
\int_0^\pi f(2\cos\theta)\,w(\theta)\,d\theta\end{equation}
 is the unique measure such that 
the polynomials \eqref{mup} are orthogonal with norm \eqref{cn}. 
By \eqref{kif},  Theorem \ref{pkt} is then equivalent to the integral formulas
\begin{subequations}\label{it}
\begin{multline}P_n(x_1,\dots,x_m)
=\frac{1}{ l!}\prod_{i=1}^n\frac{(-q;q)_i}{(q;q)_i}\prod_{i=1}^m x_i^l \prod_{1\leq i<j\leq m}(x_j-x_i)\\
\times \int_{0\leq \theta_1,\dots,\theta_{l}\leq\pi} \prod_{{1\leq j\leq m,\, 1\leq k\leq l}}(x_j+x_j^{-1}+e^{2i\theta_k}+e^{-2i\theta_k})\\
\times\prod_{1\leq j<k\leq l}(e^{2i\theta_k}+e^{-2i\theta_k}-e^{2i\theta_j}-e^{-2i\theta_j})^2\prod_{k=1}^{l} w(\theta_k)\,d\theta_k,
\end{multline}
when $n-m=2l-1$ and
\begin{multline}P_n(x_1,\dots,x_m)
=\frac{1}{ l!}\prod_{i=1}^n\frac{(-q;q)_i}{(q;q)_i}\prod_{i=1}^m x_i^l(1-x_i)
\prod_{1\leq i<j\leq m}(x_j-x_i)
\\
\times
\int_{0\leq \theta_1,\dots,\theta_{l}\leq\pi} \prod_{{1\leq j\leq m,\, 1\leq k\leq l}}(x_j+x_j^{-1}+e^{2i\theta_k}+e^{-2i\theta_k})\\
\times\prod_{1\leq j<k\leq l}(e^{2i\theta_k}+e^{-2i\theta_k}-e^{2i\theta_j}-e^{-2i\theta_j})^2\prod_{k=1}^{l} (e^{i\theta_k}+e^{-i\theta_k})^2 w(\theta_k)\,d\theta_k,
\end{multline}
\end{subequations}
when $n-m=2l$.


To prove this we will
 use the following alternative expression for $Q_\lambda$. An equivalent identity is due to Nimmo \cite[Eq.\ (A12)]{n}, see also \cite[Ex.\ III.8.13]{m}.

\begin{lemma}\label{nl}
Let $l$ be the  integral part of  $(n-m)/2$. Then,
\begin{multline*}Q_{(\lambda_1,\dots,\lambda_m)}(x_1,\dots,x_n)\\
=\frac{2^{m-l}}{l!}\prod_{1\leq i<j\leq n}\frac{x_i+x_j}{x_i-x_j}
\sum_{\sigma\in S_n}\sgn(\sigma)\prod_{i=1}^m x_{\sigma(i)}^{\lambda_i}\prod_{i=1}^l\frac{x_{\sigma(m+2i-1)}-x_{\sigma(m+2i)}}{x_{\sigma(m+2i-1)}+x_{\sigma(m+2i)}}.
\end{multline*}
In particular, if $l$ is the integral part of $(n+1-m)/2$, then
\begin{multline}\label{np}P_n(x_1,\dots,x_m)\\
=\frac{1}{2^ll!}\prod_{j=1}^n\frac{(-q;q)_j}{(q;q)_j}
\sum_{\sigma\in S_{n+1}}\sgn(\sigma)\prod_{i=1}^m x_i^{\sigma(i)-1}\prod_{i=1}^l\frac{1-q^{\sigma(m+2i)-\sigma(m+2i-1)}}{1+q^{\sigma(m+2i)-\sigma(m+2i-1)}}.
\end{multline}
\end{lemma}

\begin{proof}
By definition,
$$Q_\lambda(x)=\frac{2^m}{(n-m)!}\sum_{\sigma\in S_n}\prod_{i=1}^m x_{\sigma(i)}^{\lambda_i}\prod_{\substack{1\leq i\leq m\\1\leq i<j\leq
  n}}\frac{x_{\sigma(i)}+x_{\sigma(j)}}{x_{\sigma(i)}-x_{\sigma(j)}}. $$
We write
$$\prod_{\substack{1\leq i\leq m\\1\leq i<j\leq
  n}}\frac{x_{\sigma(i)}+x_{\sigma(j)}}{x_{\sigma(i)}-x_{\sigma(j)}}
=AB,$$
where 
$$A=\prod_{1\leq i<j\leq n}\frac{x_{\sigma(i)}+x_{\sigma(j)}}{x_{\sigma(i)}-x_{\sigma(j)}}=\sgn(\sigma)\prod_{1\leq i<j\leq n}\frac{x_i+x_j}{x_i-x_j},$$
$$B=\prod_{m+1\leq i<j\leq n}\frac{x_{\sigma(i)}-x_{\sigma(j)}}{x_{\sigma(i)}+x_{\sigma(j)}}.$$

Next, we note that
$$\prod_{1\leq i<j\leq n-m}\frac{x_i-x_j}{x_i+x_j}=\frac{1}{2^ll!}\sum_{\sigma\in S_{n-m}}\sgn(\sigma)\prod_{i=1}^l\frac{x_{\sigma(2i-1)}-x_{\sigma(2i)}}{x_{\sigma(2i-1)}-x_{\sigma(2i)}}. $$
Indeed, if $n-m$ is even this is Schur's pfaffian evaluation \eqref{spa}, while if $n-m$ is odd it is the case $x_{2m}=1$ of \eqref{spa}. It follows that the factor $B$ can be traded for
$$\frac{(n-m)!}{2^ll!}\prod_{i=1}^l\frac{x_{\sigma(m+2i-1)}-x_{\sigma(m+2i)}}{x_{\sigma(m+2i-1)}+x_{\sigma(m+2i)}}; $$
since the result is  antisymmetrized over  $S_n$ we need not perform the  antisymmetrization over the smaller group $S_{n-m}$. This completes the proof.
\end{proof}

It will be convenient to introduce the functional
$$\int f(t)\,d\lambda(t)=
\int_0^\pi \frac{f(-e^{2i\theta})-f(-e^{-2i\theta})}{e^{-2i\theta}-e^{2i\theta}}\,w(\theta)\,
d\theta,
 $$
defined on Laurent polynomials.

\begin{lemma}\label{lml}
One has
$$\int t^k\,d\lambda(t)=\frac{1-q^k}{1+q^k},\qquad k\in\mathbb Z. $$
\end{lemma}

\begin{proof}
Since 
$$\int f(t)\,d\lambda(t)=-\int f(t^{-1})\,d\lambda(t),$$
we may assume that $k$ is positive. We may then write
\begin{equation*}\begin{split}\int t^k\,d\lambda(t)&=\frac {(-1)^{k+1}}2
\int_{-\pi}^\pi \frac{e^{2i\theta k}-e^{-2i\theta k}}{e^{2i\theta}-e^{-2i\theta}}\,w(\theta)\,d\theta\\
&=\frac{(-1)^{k+1}}2 \sum_{j=0}^{k-1}
\int_{-\pi}^\pi e^{2i\theta(2j+1-k)}\,w(\theta)\,d\theta.
\end{split}\end{equation*}
By \cite[Eq.\ (7.4.18)]{gr},
$$\int_{-\pi}^\pi e^{2ik\theta}\,w(\theta)\,d\theta
=2(-1)^kq^{k}\frac{q^{-1}-q}{(1+q^{k-1})(1+q^{k+1})}=2(-1)^k(a_{k-1}-a_{k+1}),$$
where $a_k=q^k/(1+q^k)$.
Thus, the above sum telescopes, giving
$$\int t^k\,d\lambda(t)=a_{-k}-a_k,$$
which simplifies to the desired expression.
\end{proof}

We now use Lemma \ref{lml} to represent the numbers 
$$\frac{1-q^{\sigma(m+2i)-\sigma(m+2i-1)}}{1+q^{\sigma(m+2i)-\sigma(m+2i-1)}}$$
 appearing in \eqref{np} as integrals. This gives
\begin{multline*}P_n(x_1,\dots,x_m)\\
=\frac{(-1)^l}{2^l l!}\prod_{j=1}^n\frac{(-q;q)_j}{(q;q)_j}
\sum_{\sigma\in S_{n+1}}\sgn(\sigma)\prod_{i=1}^m x_i^{\sigma(i)-1}\prod_{i=1}^l
\int (-t)^{\sigma(m+2i)-\sigma(m+2i-1)}\,d\lambda(t).
\end{multline*}

Interchanging summation and integration, the sum becomes a Vandermonde determinant, for which we introduce the notation
$$\Delta(x)=\det(x_i^{j-1})=\prod_{i<j}(x_j-x_i). $$
If $n-m=2l-1$, we thus obtain
\begin{multline*}P_n(x_1,\dots,x_m)
=\frac{1}{2^l l!}\prod_{j=1}^n\frac{(-q;q)_j}{(q;q)_j}\\
\times\int \Delta(x_1,\dots,x_m,t_1^{-1},t_1,\dots,t_l^{-1},t_l)\,d\lambda(t_1)\dotsm d\lambda(t_l),
\end{multline*}
 whereas if $n-m=2l$, 
\begin{multline*}P_n(x_1,\dots,x_m)
=\frac{1}{2^l l!}\prod_{j=1}^n\frac{(-q;q)_j}{(q;q)_j}\\
\times\int \Delta(x_1,\dots,x_m,t_1^{-1},t_1,\dots,t_l^{-1},t_l,1)\,d\lambda(t_1)\dotsm d\lambda(t_l).
\end{multline*}

Next, we  apply the following elementary identities.

\begin{lemma}
One has
\begin{multline*}\Delta(x_1,\dots,x_m,t_1,t_1^{-1},\dots,t_l,t_l^{-1})\\
=\prod_{i=1}^m x_i^l
\prod_{j=1}^l(t_j-t_j^{-1})\,\Delta(x)
\Delta(\tau)^2\prod_{1\leq j\leq m,1\leq k\leq l}(\xi_i-\tau_j),
 \end{multline*}
\begin{multline*}\Delta(x_1,\dots,x_m,t_1,t_1^{-1},\dots,t_l,t_l^{-1},1)\\
=\prod_{i=1}^m x_i^l(1-x_i)
\prod_{j=1}^l(t_j-t_j^{-1})(2-\tau_j)\,\Delta(x)
\Delta(\tau)^2\prod_{1\leq j\leq m,1\leq k\leq l}(\xi_i-\tau_j),
 \end{multline*}
where $\xi_i=x_i+x_i^{-1}$, $\tau_i=t_i+t_i^{-1}$.
\end{lemma}

 Writing
$$\int (t_j-t_j^{-1})f(\tau_j)\,d\lambda(t_j)=2\int_0^\pi
f(-e^{2i\theta_j}-e^{-2i\theta_j})\,w(\theta_j)\,d\theta_j,
 $$
we  recover  \eqref{it}.
This completes our alternative proof of Theorem \ref{pkt}.

\end{document}